%% file: opti_nonclassical_TAC_twocol.tex
\documentclass[journal]{IEEEtran}  

\IEEEoverridecommandlockouts                              
\usepackage{multicol}
\usepackage{wrapfig}
\input{macros-ieee.tex}
\usepackage{multirow}

\title{\LARGE \bf 
An Optimizer's Approach to Stochastic Control Problems with Nonclassical Information Structures
}

\author{Ankur A. Kulkarni$^{1}$ and Todd P. Coleman$^{2}$
\thanks{$^{1}$Ankur is with the Systems and Control Engineering group at the Indian Institute of Technology Bombay, Mumbai, India 400076.
        {\tt\small kulkarni.ankur@iitb.ac.in}. The contents of this paper were presented in part at the IEEE Conference on Decision and Control, 2012 at Maui, Hawaii, USA \cite{kulkarni12optimizer}. }%
\thanks{$^{2}$Coleman is with the Department of Bioengineering, University of California, San Diego, La Jolla, CA USA.
        {\tt\small tpcoleman@ucsd.edu}.}%
}

\begin{document}

\maketitle
\thispagestyle{empty}
\pagestyle{empty}

\begin{abstract}
We present an optimization-based approach to stochastic control problems with nonclassical information structures. We cast these problems equivalently as optimization problems on joint distributions. The resulting problems are necessarily \textit{nonconvex}. Our approach to solving them is  through \textit{convex relaxation}. We solve the instance solved by Bansal and Ba\c{s}ar~\cite{bansal87stochastic} with a particular application of this approach that uses the data processing inequality for constructing the convex relaxation. Using certain $f$-divergences, we obtain a new, larger set of inverse optimal cost functions for such problems.  Insights are obtained on the relation between the structure of cost functions and of convex relaxations for inverse optimal control. 
\end{abstract}
\def\td{\textrm{-}}
\def\BB{Bansal-Ba\c{s}ar\xspace}
\def\Qscrfdpi{\Qscr_{f\td{\rm DPI}}}

\section{MOTIVATION AND CONTRIBUTION}\label{sec:motivation}
\begin{figure}[h]
 {\centering
 \includegraphics[scale=.33]{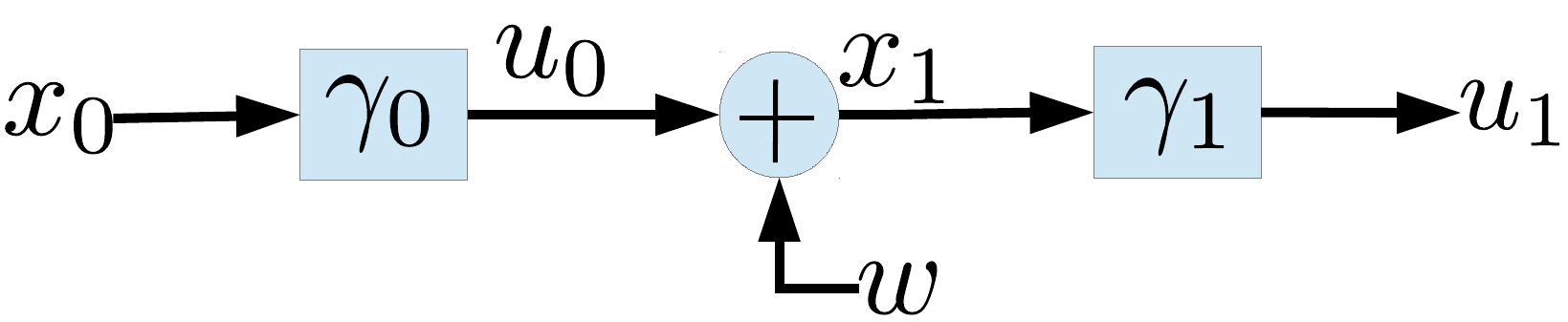}
 \caption{Nonclassical information structure and the setting of $\Nbf$}
 \label{fig:nonclassical} }
 \end{figure}

This paper is motivated by the following stochastic control problem with nonclassical information structure:
  minimize 
\begin{align}
J(\gamma_0,\gamma_1) &= \Ebb \left[ \kappa(x_0,u_0,x_1,u_1) \right], \tag{$\Nbf$} \label{eq:n} \\
u_0 &= \gamma_0(x_0), \quad x_1 = u_0 + w, \quad u_1 = \gamma_1(x_1), \non
\end{align}
where $x_0,w$ are random variables with fixed statistics, 
$u_0,u_1$ are constrained by the \textit{information structure}, \ie, $u_0$ is adapted to $x_0$ alone and $u_1$ is adapted to $x_1$ alone,   
and the decision variables are measurable functions $\gamma_0,\gamma_1$ (Fig \ref{fig:nonclassical}). 
In general, these problems are challenging to solve.
For example, if 
$\kappa(x_0,u_0,x_1,u_1) = (x_0-u_0)^2 + (u_0-u_1)^2,$ 
and $x_0,w$ are Gaussian, then this becomes the long standing open problem of Witsenhausen~\cite{witsenhausen_counterexample_1968} about which many fundamental issues remain to be clarified.

In 1987, Bansal and Ba\c{s}ar~\cite{bansal87stochastic} solved a nontrivial problem with this information structure. They consider the case where the problem is to minimize
\begin{align}
  \kappa(x_0,u_0,x_1,u_1) &= k_0 u_0^2 + s_{01} u_0x_0 + \left( u_1-x_0 \right)^2, \tag{$\Bbf$}\label{eq:bbj} 
\end{align} where $x_0 \sim \Nscr(0,\sigma_0^2), w \sim \Nscr(0,\sigma_w^2)$ are independent Gaussian random variables and $k_0>0$ is a scalar constant.  Their approach was to interpret Fig \ref{fig:nonclassical} as a communication system where $\gamma_0$ and $\gamma_1$ are the encoder and decoder, respectively, and then use information-theoretic bounds to obtain a solution.
An important aspect of their proof is that it relies critically on the information-theoretic data processing inequality (DPI). This remains the only known logic to proving optimality for problem $ \Bbf $~\cite{basar08variations}.

Upon studying the proof of Bansal and Ba\c{s}ar, it is not at once clear if their proof is a result of a mathematical coincidence, or if it can be generalized in a way whereby it can be applied to other problems with this information structure. 
In this paper we would like to \textbf{a)} understand if there is some aspect of this proof that can be generalized to arrive at a methodology for other instances of $ \Nbf $.  It is evident that this would involve understanding the role of the DPI in this proof. 
The similarities between $ \Bbf $ and the Gaussian test channel suggest that the DPI arises from communication-theoretic considerations. Therefore,   
  we would like to \textbf{b)} understand if the role of the DPI in the solution of $ \Bbf $ is primarily communication-theoretic or if there is another interpretation that makes the role more general.
  And \textbf{c)} we would like to know if there is an explanation of what makes $ \Bbf $, which is prima facie only a slight modification of the Witsenhausen problem~\cite{bansal87stochastic},  tractable, even while the Witsenhausen problem itself remains unsolved.

This paper argues that the solution of $ \Bbf $ can be obtained as a special case of a more general optimization-based paradigm and the role of the DPI is precisely that of convexifying a nonconvex problem.  We cast $\Nbf$  as an equivalent problem of minimization of the cost over the joint distribution of all variables, but  where $\gamma_0$ and $\gamma_1$ (the encoder and decoder) are allowed to be {\em random}. The resulting problem has a linear objective and necessarily \textit{nonconvex} constraints. 
Nonconvexity implies there is no general approach for solving $\Nbf$. A possible solution strategy is to construct a \textit{convex relaxation} with the same objective and a larger and convex feasible region. This problem being convex is potentially more accessible. The strategy then is to find a solution of the relaxation that is feasible for $\Nbf$, which would thereby be a solution of $\Nbf$.

\subsection{Main new results and insights}
Our main finding is that $\Bbf$ is solvable by the above approach. We show that $ \Bbf $ is a nonconvex problem, but there is a convex relaxation of $\Bbf$ which can be solved to obtain a solution that is feasible for $\Bbf$. The construction of the relaxation involves the DPI. 

Implicit in the proof of Bansal and Ba\c{s}ar is the use of the variational equations for the rate-distortion and capacity-cost function and their relation using the DPI. 
While we also use the DPI, our proof shows that \textit{only its convexity properties are needed} in solving $\Bbf$; no communication-theoretic properties are required. Consequently, the DPI is only an artifice for convexifying the problem and as such may be replaced by any other suitable means of convexification. 
Since the DPI is only \textit{one} candidate convexification,  it is plausible that other convexifications would yield other solutions of $ \Bbf $ or solutions to other variants of $ \Nbf $. 

Indeed such solutions are readily obtainable for the problem of \textit{inverse optimal control} in which \textit{cost functions}, for which a candidate policy is optimal, are sought. This problem is useful when one knows a policy adopted by  controllers (suggested by, say, experimental data), and these controllers are known to be acting ``optimally'' according to some unknown criterion which one would like to determine (see \eg,~\cite{abbeel2010autonomous}). 
With the convex relaxation viewpoint,  we show that inverse optimal control requires the perfect agreement of the cost function and the function that convexifies the problem. The standard set of inverse optimal cost functions obtained via information-theoretic arguments correspond to those with the convexifcation using the DPI. However using other divergences, such as certain $f$-divergences, we obtain a new, larger set of inverse optimal cost functions. As such, these results provide a general purpose framework to propose other convexifications that will lead to similar formalisms.

 Convex relaxation is an alternative to information theory for obtaining lower bounds for problems like $ \Nbf $. To quote Ho \etal~\cite[p.\ 307, footnote 5]{ho_teams_1978} \textit{``There exist no general sufficiency conditions to verify the optimality of a solution (of a problem like $ \Nbf $) besides the Shannon bounds''}. The information-theoretic way of obtaining a lower bound involves embedding $ \Nbf $ into a problem that permits infinite delay (\ie, arbitrary block lengths). This is somewhat unsatisfactory for two reasons: first, the appearance of mutual information is a likely consequence of the infinite block length problem, whereby the method is somewhat inflexible. And, secondly, the bounds will often be too loose (see Witsenhausen's discussion on this topic~\cite[Section 8]{witsenhausen1975intrinsic} and Ho \etal~\cite{ho_teams_1978}). Convex relaxation provides a broader paradigm for obtaining bounds for $ \Nbf $. And for problem $ \Bbf $,  
which is the most general problem we know of to have been solved by the above information-theoretic approach,
the convex relaxation based bounds include the  information-theoretic bounds as a special case.


Specifically for problem $ \Bbf, $ our results amount to a new proof of optimality  -- our proof does not use the Cauchy-Schwartz inequality or parametrization (which are two of the steps in \cite{bansal87stochastic}). The tractability of $ \Bbf $ is explained by the fact that even though $ \Bbf $ is equivalent to a nonconvex optimization problem, it admits a convex relaxation that is tight for it. Furthermore, we show that if a certain nondegeneracy condition holds, equality in the DPI is a necessary condition for optimality  in $\Bbf$. 
Consequently, any other proof of optimality of $\Bbf$ must imply this equality, even though the DPI may not have been explicitly employed in the proof.  This partially settles the open problem~\cite{basar08variations} of ascertaining whether there is a way of solving $\Bbf$ that does not use the DPI. 

Finally, by  convex relaxation 
  we give a lower bound on the Witsenhausen problem and 
 a result that indicates that a certain ``easy'' solution approach cannot succeed for it.


\subsection{Theoretical explanation of the results} \label{sec:message}
There exists an optimization-based approach for Markov decision processes (stochastic control problems with {\em classical} information structures), where one optimizes the expected cost over  distributions on state-action spaces, called \textit{occupation measures} (see \eg,~\cite{bertsekas05dynamic,borkar_topics_1991}). In this formulation, the constraints on an occupation measure are that a) it be a valid distribution and b) that it be an invariant distribution. The relaxation  of deterministic codes to random codes, is a logical extension of this approach to  problem $\Nbf$. We minimize the expected cost over distributions that are consistent with the given marginals and that respect the information structure of the problem. But while the constraints in the classical information structure result in a 
convex (in fact, linear) program, the nonclassicality of the information structure implies the nonconvexity of the corresponding optimization problem\footnote{One may think of this as a different perspective on the observation that problems with nonclassical information structure are harder than those with classical information structure.}. 

Problem $\Bbf$ of Bansal and Ba\c{s}ar is also nonconvex, but it happens to share its solution with a specific relaxation. This is because there is a convex set that contains the feasible region of $\Bbf$ and whose normal at the solution of $\Bbf$ ``agrees'' with the cost of $\Bbf$. To explain this, suppose a solution of the nonconvex problem \eqref{p:nonconvex} where $v$ is a vector and $C'$ is a nonconvex set, is also a solution 
of its convex relaxation \eqref{p:convex} where $C$ is a convex set containing $C'$.
\noindent\begin{tabular}{p{.22\textwidth}p{.22\textwidth}}
{\begin{align} %
\minimize{z \in C'} v^\top z,    \label{p:nonconvex}
\end{align}} 
& 
{\begin{align} %
\minimize{z \in C} v^\top z, \label{p:convex}%
\end{align}}\vspace{-5pt}
\end{tabular}

\begin{figure} \centering
  \vspace{-15pt}
  \includegraphics[scale=.27]{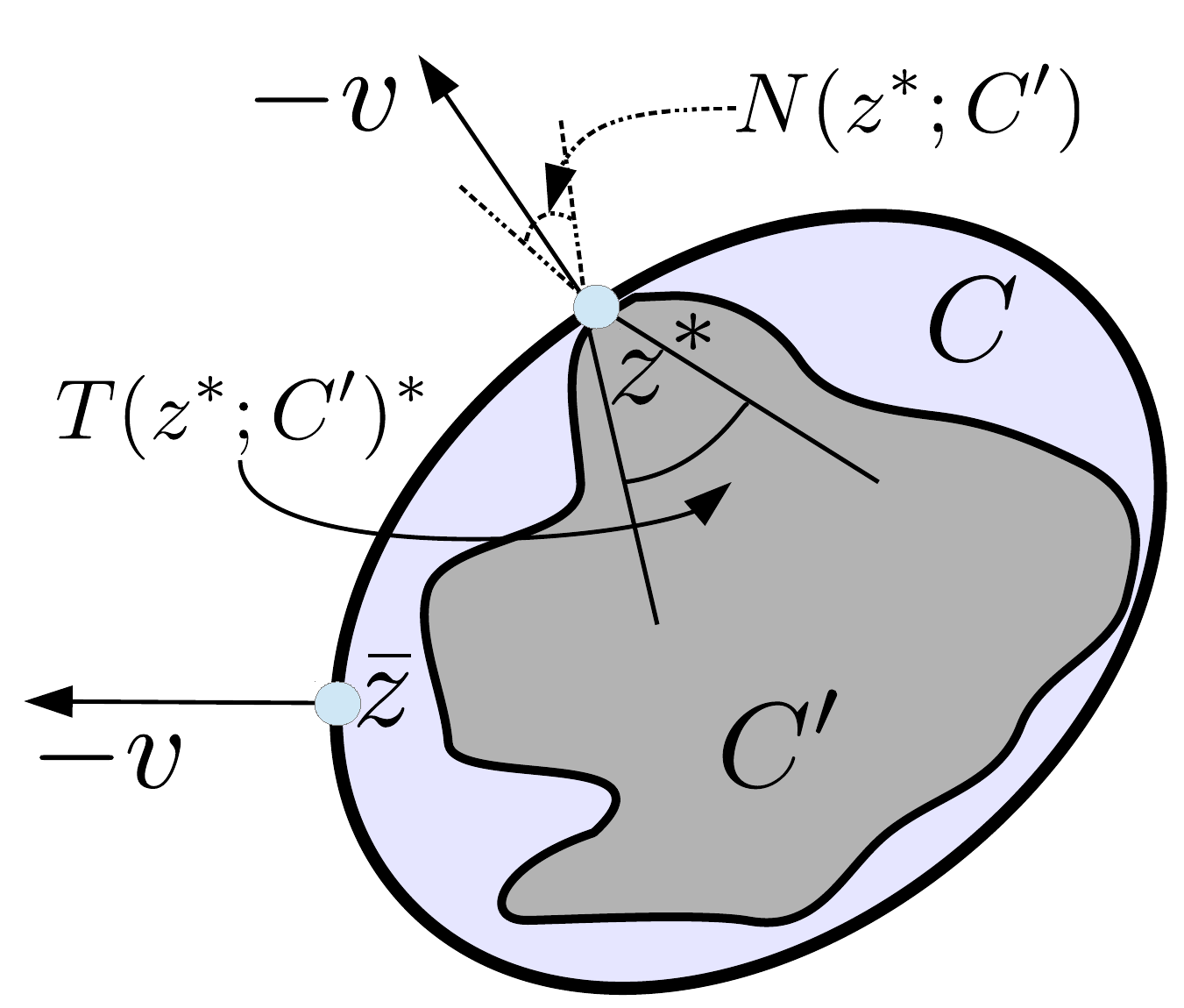}
  \caption{The point $z^*$ is optimal for both \eqref{p:convex} and \eqref{p:nonconvex}. $\bar{z}$ is optimal for \eqref{p:convex} but is not feasible for \eqref{p:nonconvex}.}
  \label{fig:optimality} 
  \vspace{-10pt}
  \end{figure}

\noindent A point $z^* \in C$ (see Fig \ref{fig:optimality}) is a solution of \eqref{p:convex} if and only if $-v$ is normal to $C$ and pointing outwards (technically, $-v$ belongs to the \textit{normal cone} of $C$ at $z^*$; see~\cite{rockafellar09variational}).   For \eqref{p:nonconvex} and \eqref{p:convex} to share a solution, the point where $-v$ is normal to $C$ must lie in $C'$. This coincidence occurs in problem $\Bbf$ if the set $C$ is formed using the DPI.  In general, an optimal solution of $C$ may not lie in $C'$ (as depicted by point $\bar{z}$ in Fig \ref{fig:optimality}). 

The notion of \textit{inverse optimality} seeks cost functions (the vector $v$) for which a candidate solution is optimal. For \eqref{p:nonconvex}, one can in general give only \textit{necessary} conditions for inverse optimality. 
KKT conditions yield a set of functions (denoted by $T(z^*;C')^*$, the \textit{dual of the tangent cone}, in Fig \ref{fig:optimality}) that is larger than the set of inverse optimal functions. 
For $v$ to be inverse optimal for \eqref{p:nonconvex}, it is sufficient that $-v$ lies in $N(z^*;C')$ (the normal cone of  $C'$ at $z^*$, cf. Fig \ref{fig:optimality}). Since $C \supset C'$, at $z^*$,  we have the relation
\begin{equation}
N(z^*;C)\subseteq N(z^*;C'). \label{eq:normalcone}
\end{equation} 
If $z^*$ solves \eqref{p:convex}, cost functions (which are now inverse optimal for \eqref{p:convex}) are guaranteed to be  inverse optimal for \eqref{p:nonconvex}. In general these form a strict \textit{subset} of the set of functions that are inverse optimal for \eqref{p:nonconvex}. Notice that \eqref{eq:normalcone} holds for any convex set $C$ containing $C'$. As such the normal cone of any convex set $ C $ provides a subset of the normal cone of $ C' $. Thus one gets a suite of inverse optimal cost functions by varying the convex set employed for creating the relaxation. 

In our case, the convexification with the DPI provides one set of inverse   optimal cost functions 
 for a variant of $\Bbf$ (with $s_{01}=0$). 
In fact it provides a stronger form of inverse optimality where `distortion-like' and `cost-like' terms are separately
 obtained. Suppose one asks for $v_1$ and $v_2$ where 
 $v_1,v_2$ lie in given orthogonal subspaces such that $v = v_1 + v_2$ is 
  inverse optimal for \eqref{p:nonconvex}. These are obtained as the projection of the normal cone along these subspaces, which can be obtained cleanly if $C$ has a favorable structure. In $\Bbf$ with $s_{01}=0$, the DPI allows the inverse optimal distortion and cost to be separately unravelled. This is because of the agreement between the structure of the DPI and the structure of these terms.

By replacing the KL divergence with other suitable $ f $-divergences, one can vary the choice of the set $ C$  and thereby obtain a larger set of inverse optimal cost functions, wherein $ f $ is also a parameter defining these functions. 
\subsection{Background and organization}
There are two main contributions from our work. First, convex relaxation can be seen more generally as an \textit{approach} for problems with 
nonclassical information structure.
Second, in the particular case of problem $ \Bbf $ our work clarifies certain issues, including its hardness and the role of the DPI in its solution.    Here we present some background relevant to these aspects of our work. Problems like $ \Nbf $ have a long and rich history and it would not be possible to cover all of it in the space available. 

Problem $ \Nbf $ is a stochastic control problem where the information available to the latter-acting decision maker is dependent on the actions of the earlier-acting decision maker. The information structure of $ \Nbf $ is thus \textit{dynamic} in the terminology of Ho~\cite{ho1980team}. In addition, the latter-acting decision maker $\gamma_1$ does not know what the earlier-acting decision maker $\gamma_0$ knows, whereby the information structure of $ \Nbf $ is not \textit{partially nested}~\cite{ho1980team,basar99dynamic}. The partially nested information structure is essentially the only one that is equivalent to a static information structure~\cite{ho80another}. In that respect, $ \Nbf $ is perhaps the simplest example of a problem with a dynamic  information structure that is not partially nested.

The famous counterexample of Witsenhausen~\cite{witsenhausen_counterexample_1968}, a special case of $ \Nbf $,  showed that linearity of state evolution, quadratic cost and Gaussian noise, do not ensure a linear optimal controller if the information structure is not classical (for the same cost function, if $ u_1 $ is allowed to depend on $ x_0,u_0,u_1 $, one gets a linear optimal controller).   
This has prompted a line of research asking, ``When is an optimal controller linear?''~\cite{kalman1964linear}, and the result of Bansal and Ba\c{s}ar~\cite{bansal87stochastic} was one of the most general results along this direction for the information structure of $ \Nbf $. Prior to \cite{bansal87stochastic} it was known that for the Gaussian test channel (\ie $ \Bbf $ with $ s_{01} =0$), the optimal codes are linear.  Bansal and Ba\c{s}ar noticed that this can also be extended to the case where $ s_{01}\neq 0.$ In a later paper~\cite{bansal87solutions} they consider an extension of this problem where the channel is vector-valued. In yet another extension~\cite{bansal89simultaneous}, they consider a system with dynamics and feedback. In both problems \cite{bansal87solutions,bansal89simultaneous}, the optimal controllers are found to be linear. In each of these proofs the DPI was 
used as a key ingredient. It remains an open problem to ascertain if there is another way of proving optimality for $ \Bbf $ that does not use the DPI~\cite{basar08variations}. In this context, our work shows that in the solution of $ \Bbf $ the DPI is only an artifice for convexifying the problem, and as such may be replaced by any other suitable convex set. However equality in the DPI is a necessary condition for optimality in $ \Bbf, $
so any other proof must imply this equality.

In problems like $ \Nbf,$ a fundamental challenge is of understanding what makes these problems intractable. An attempt in this direction can be made by observing that since the argument of $ \gamma_1 $ is dependent on the value of $\gamma_0$, and $ \gamma_1 $ does not have access to the argument of $ \gamma_0 $, problem $ \Nbf $ is essentially an optimization of  $J(\gamma_0,\gamma_1(\gamma_0))$ over $ \gamma_0,\gamma_1$~\cite{ho1980team}. This is clearly a nonstandard problem in optimization over functions. Furthermore, it is a \textit{nonconvex} problem in $ \gamma_0 $ whenever $ \gamma_1 $ is not linear or a suitable function~\cite{ho1980team,witsenhausen_counterexample_1968}. 
To quote Ho~\cite[p.\ 648]{ho1980team}, \textit{``The computational as well as theoretical difficulties sketched above makes (P3) (\ie, $ \Nbf $) a most difficult problem. It is still unsolved today, some 12 years\footnote{Now, 45 years ...} after Witsenhausen first analyzed it. In fact, ... we do not even have a sufficient condition for optimality for (P3) similar to the usual second variation condition in function optimization. Optimization problem of the type $ J(\gamma_1,\gamma_2(\gamma_1)) $ involving composition of the optimizing functions simply have not been studied with any kind of systematic effort.''
}\footnote{On another note, we recall the effort of Papadimitriou and Tsitsiklis~\cite{papadimitriou85intractable}  to understand the hardness of the Witsenhausen problem through the lens of computational complexity.}


The above explanation for the hardness of $ \Nbf $ is somewhat unsatisfactory since it does not 
explain the differences in the degree of hardness of nonclassical problems. Specifically, note that while we know little about the solution of the Witsenhausen problem, for the Bansal-Ba\c{s}ar problem $ \Bbf $, which is a slight variant of the Witsenhausen problem, we do know the solution. This raises the natural question: what makes the Bansal-Ba\c{s}ar problem tractable? To the best of our knowledge, an answer to this question from the viewpoint of optimization over functions is not known. Bansal and Ba\c{s}ar~\cite{bansal87stochastic} have made the observation that the presence of a term in `$ u_0 u_1 $' makes the Witsenhausen problem hard. However this does not yet explain what makes $ \Bbf $ tractable. 

Our work shows that $ \Nbf $ is equivalent to a nonconvex optimization problem in the space of joint distributions. Importantly, this result does not rely on the structure of the cost $ \kappa $, but is instead a consequence of the information structure of $ \Nbf $. The tractability of $ \Bbf $ is then explained by the fact that, although it is a nonconvex problem, there happens to be a convex relaxation that is tight for it.

A closely related paradigm where one optimizes over joint distributions that are constrained to be consistent with given marginals is the subject of optimal transport theory~\cite{villani2008optimal}. The framework presented here may be seen as optimal transport with \textit{additional constraints} which arise from the information structure. 
Recent work of Wu and Verd\'{u}~\cite{wu11optimal} also casts the problem as an optimization over distributions, but the eventual problem they address has a different formulation from ours. Their focus is specifically on the Witsenhausen problem; they exploit the quadratic structure of the cost function and their newly discovered properties of the MMSE estimator~\cite{wu2012functional}.  With this, they can relate the Witsenhausen problem to a standard optimal transportation problem with quadratic cost, to  characterize properties of an optimal solution.

On the other hand, our contributions are centered around the information
structure, the nonconvexity of the problem and addressing it through convex relaxation.

At the heart of problems like $ \Nbf $ is the interplay of information and control~\cite{witsenhausen_separation_1971}, signalling and the dual role of control~\cite{basar08variations}. We are not explicitly concerned with these aspects of $ \Nbf $. However we take note of the fact that $ \Nbf $ has inspired a search for an understanding of the role of information in control~\cite{mitter1999information,grover2010witsenhausen}. Indeed, a related line of research attempts to do dynamic programming for problems with various information structures~\cite{mahajan2009optimal}. For other perspectives on this problem, we refer the reader to the recent tutorial~\cite{mahajan2012information}.

This paper is organized as follows. In Section \ref{sec:bb}, we recall the proof of 
Bansal and Ba\c{s}ar and the communication-theoretic approach to problems like $ \Nbf $. Section \ref{sec:opti} contains the optimization formulation, convex relaxation and the solution of $\Bbf$. Section  \ref{sec:ioc} concerns inverse optimal control and Section \ref{sec:witsen} concerns the convex relaxation of the Witsenhausen problem. We conclude in Section \ref{sec:conclusions}.

\section{THE PROBLEM OF BANSAL AND BA\c{S}AR }\label{sec:bb}
In this section we recall the proof of Bansal and Ba\c{s}ar. In problem $\Bbf$, $x_0 \sim \Nscr(0,\sigma_0^2)$ and $w \sim \Nscr(0,\sigma_w^2)$ (cf. Fig \ref{fig:nonclassical}) are independent random variables. 
$\Bbf$ is essentially the Witsenhausen problem but with the coefficient of $u_1u_0$ taken to be zero~\cite{bansal87stochastic}. 
Bansal and Ba\c{s}ar interpret this problem as that of minimizing the distortion $(u_1-x_0)^2$ in a communication system subject to soft constraints on power $u_0^2$ and $u_0x_0$. Their proof \cite[section III]{bansal87stochastic} uses the following series of arguments: 
\begin{enumerate}
\item[\textbf{1.}]  The DPI implies $I(x_0;u_1) \leq I(u_0,u_0+w)$.
\item[\textbf{2.}] Standard results about Gaussian channels give that, under the constraint $\Ebb[u_0^2] \leq P^2$, 
\begin{align}
\frac{1}{2}\log \frac{\sigma_0^2}{\Ebb[(u_1-x_0)^2]} &\leq   I(x_0;u_1)\label{eq:bbbound1} \\
\leq I(u_0,u_0+w) &\leq \frac{1}{2}\log \frac{P^2 + \sigma_w^2}{\sigma_w^2}. \label{eq:bbbound2}
\end{align} 
Combining these gives a lower bound on the distortion: $\Ebb[(u_1-x_0)^2] \geq \frac{\sigma_0^2\sigma_w^2}{P^2 + \sigma_w^2}.$
\item[\textbf{3.}] Cauchy-Schwartz inequality gives a lower bound on the third term
$    s_{01}\Ebb\left[ u_0x_0\right]  \geq  -|s_{01}|P\sigma_0.$
\end{enumerate}
Thus the optimal value of $ \Bbf $, $J^*$ is bounded below as $J^* \geq   \frac{\sigma_0^2 \sigma_w^2}{{P^*}^2 + \sigma_w^2} + k_0 {P^*}^2 - |s_{01}|P^*\sigma_0,$ where $P^*$ satisfies
\begin{equation}
 (2k_0P^* - |s_{01}|\sigma_0)({P^*}^2 + \sigma_w^2)^2 = 2P^*\sigma_0^2\sigma_w^2. \label{eq:bbpstar}
\end{equation} 
Finally, this bound is shown to be tight by the explicit construction of maps $\gamma_0^*,\gamma_1^*$ for which 
$J(\gamma^*_0(x_0),\gamma^*_1(x_1))$ equals the lower bound. The most important step here is obtaining the lower bound on the distortion. The only known approach for this is through the DPI. 
\subsection{Communication theoretic explanation and information theoretic bounds}
\begin{figure}[h]\vspace{-5pt}
\centering
\includegraphics[scale=.25]{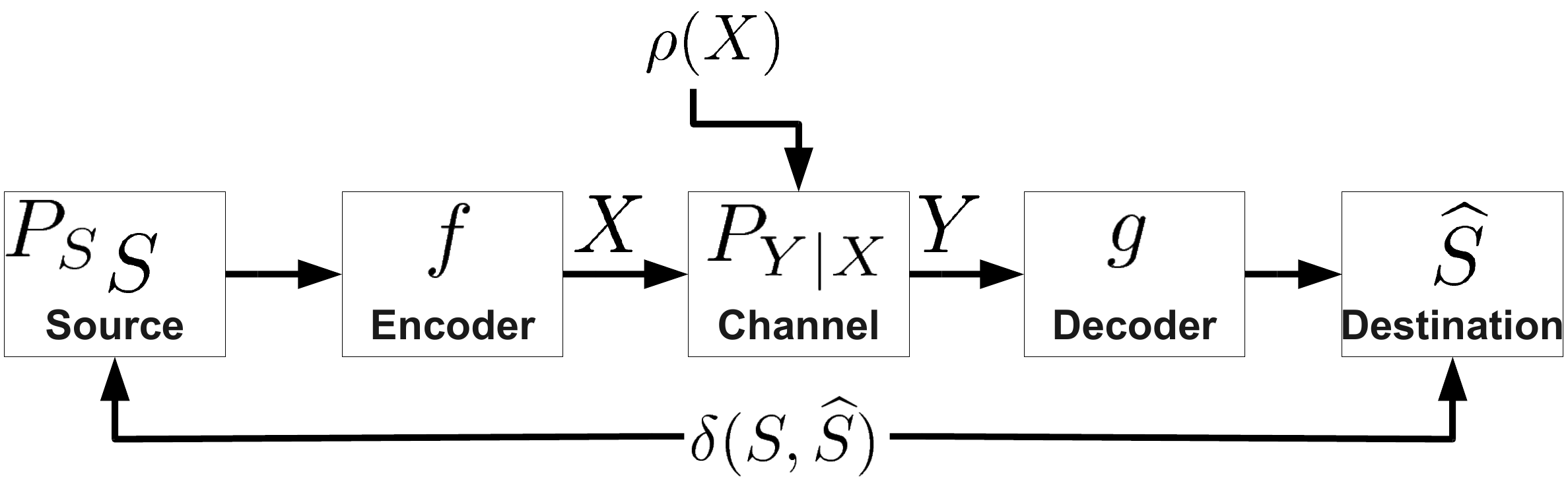}
\caption{Communication system with performance metrics}
\label{fig:comm}
 \vspace{-10pt}
\end{figure}

We now  explain the proof of Bansal and Ba\c{s}ar by comparing problem $\Bbf$ with a problem from communication theory. 
We use the notation of Fig \ref{fig:comm}. 
A source generates symbols (random variables) $S$ taking values in a space $\Sscr$, with probability $P_S(S)$. An encoder $f$ 
maps these symbols to a space $\Xscr$. Each symbol $X \in \Xscr$ passes through a channel, which produces a symbol $Y \in \Yscr$ with
probability $P_{Y|X}(Y|X)$. A decoder $g$ maps $Y$ to a
destination space $\Sscrhat$. The system is endowed with two performance metrics: a {\em distortion measure} $\delta: \Sscr \times \widehat{\Sscr} \rightarrow [0,\infty)$ and a {\em cost} $\rho: \Xscr \rightarrow [0,\infty)$. 

Standard communication theory employs \textit{block codes} wherein multiple symbols are 
coded together. An $n$-block encoder (resp., decoder) is a map $f: \Sscr^n \rightarrow \Xscr^n$ (resp., $g:\Yscr^n \rightarrow \Sscrhat^n$). A \textit{single-letter code} is a 1-block code. The subject of rate-distortion theory is the minimization 
average distortion and average cost, over the class of codes of arbitrary block length. Only in certain exceptional cases does it turn out that the optimal code is a single-letter code. Instead of seeking such cases, one may ask the inverse question studied by Gastpar \etal~\cite{gastpar03tocode}: when is a single-letter code Pareto optimal for average distortion and cost, over all codes of all block lengths?

For every single-letter code one can find a class of \textit{inverse optimal} cost and distortion functions such that the code is optimal for them in the above sense. To characterize this class, Gastpar uses the communication-theoretic constructs of the {\em rate-distortion} function $R$ and {\em capacity-cost} function $C$. For each $D \geq 0$, $R(D)$ 
is given by 
\begin{align}
 R(D) = \min_{p(\shat,s)\ :\ \Ebb\left[ \delta(S,\Shat)\right] \leq D,\ p(s) = P_S(s)} I(S,\Shat). \label{eq:rd}
\end{align}
where $I(S,\Shat)$ is the mutual information between $S$ and $\Shat$. For each $P \geq 0$, the capacity-cost function 
is \begin{align}
 C(P) = \max_{p(x,y):\Ebb\left[ \rho(X)\right] \leq P,\ p(y|x) = P_{Y|X}(y|x)} I(X,Y). \label{eq:cp}
\end{align}
where $I(X,Y)$ is the mutual information between $X$ and $Y$.

If $D = \frac{1}{n}\sum_k \Ebb [\delta(S_k,\Shat_k)]$ and $P = \frac{1}{n}\sum_k \Ebb[\rho(X_k)]$ are the average distortion and cost incurred by an $n$-block code, then $R(D) \leq C(P)$. This follows from the DPI~\cite{cover06elements} 
\begin{equation}
 I(X,Y) \geq I(S,\Shat). \label{eq:dpi}
\end{equation}
Gastpar argues that, except in some degenerate cases, 
$R(\Delta\s)=C(\Gamma\s)$, is \textit{equivalent} to the Pareto optimality of average distortion and cost values $\Delta\s$, $\Gamma\s$ (cf. Lemma 1, \cite{gastpar03tocode}). 

If $R(\Delta\s)=C(\Gamma\s)$, the \textit{Pareto optimal single-letter code} comprises of maps $f,g$ such that
$f(S) \sim p^*(x)$ and $g(Y)|S \sim p^*(\shat|s)$ where $p^*(x)$ and $p\s(\shat|s)$ achieve the optima in \eqref{eq:cp},\eqref{eq:rd}, respectively. Note that the \textit{existence} of such a single-letter code is not guaranteed by this argument. 
If one posits distributions $p^*(x)$ and $p^*(\shat|s)$ induced by a single-letter code, 
inverse (Pareto) optimal cost and distortion are given by~\cite[Lemma 3, 4]{gastpar03tocode} 
\begin{align}
\rho(x)  &= c_1 D(P_{Y|X}(\cdot|x)|| p\s_Y(\cdot)) + \rho_0   + \beta(x)\I{p^*(x)>0} \non\\
 \delta(s,\shat) &= -c_2 \log p\s(s|\shat) + d_0(s) \label{eq:dchar}
\end{align}
where $c_1,c_2>0$, $\rho_0 \in\Real$, $d_0: \Sscr \rightarrow \Real$ is arbitrary, $\beta : \Xscr \rightarrow [0,\infty)$ 
and $p^*_Y$ is the marginal of $Y$ under $p^*$. We recall that this characterization is also known from~Csiszar and K{\"o}rner~\cite{csiszar2011information}.

\subsubsection{Relation to the \BB problem} \label{sec:gastparbb}
The key distinction between the information theoretic problem described above and the control problem $ \Nbf $ is that the latter assumes a fixed (effectively, unit) block length. Specifically, since the controllers $\gamma_0,\gamma_1$ in problem $\Nbf$ act on a single sample 
of their respective arguments, they are single-letter codes. One may consider an ``information theoretic'' version of $ \Nbf $ in which the code allows infinite delay (\ie arbitrary block lengths) and the objective is the average cost. The achievability results from rate-distortion theory when applied to this problem  do not in general provide for the existence of a single-letter code (recall the discussion in \cite{witsenhausen1975intrinsic} and \cite{ho_teams_1978}). 

The cost in $\Bbf$ is a sum of a distortion-like term $(u_1 - x_0)^2$, a cost-like term $k_0u_0^2$ and a generalized cost term $s_{01}u_0 x_0$. Suppose $s_{01}=0$. In this case problem $\Bbf$ asks for a single-letter code that minimizes a sum of distortion and cost \textit{over the class of all single-letter codes}. Clearly, if a single-letter code exists that is Pareto optimal in the sense of Gastpar, it is (upon appropriate scaling of cost functions) also optimal for $\Bbf$. Problem $\Bbf$ with $s_{01}=0$  happens to be the exceptional case where such a code does exist. Indeed in \eqref{eq:dchar} if $f$ is taken to be the identity and $\Shat =g(Y)= \frac{\sigma_0^2}{\sigma_0^2 + \sigma_w^2}Y,$  we get $\rho(x) = c_1 x^2 + \rho_0$ and $\delta(s,\shat)  = c_2(s -\shat)^2 +d_0(s)$~\cite{gastpar03tocode},  which matches the distortion and cost terms in $\Bbf$. 

With this understanding the appearance of the data processing  inequality  in the proof of Bansal and Ba\c{s}ar can be explained as a vestige of the communication-theoretic problem that allows infinite delay and a stronger notion of optimality  that seeks optimality of a single letter code over codes of arbitrary block lengths. 

Importantly, if, to obtain inverse optimal cost functions for the control problem $ \Nbf $, one argues via Gastpar's notion of inverse optimality, then the DPI is indespensable: a single-letter code is optimal in the sense of Gastpar \textit{if and only if} it satisfies the characterizations in \eqref{eq:dchar}.  The optimization approach we develop in the following sections is comparatively more flexible, and includes in it as a special case, the bound/inverse optimality characterizations obtained above.

\section{THE OPTIMIZATION-BASED APPROACH}\label{sec:opti}
This section contains the main contributions of this paper. In Section \ref{sec:formulation}, we formulate  the problem $\Nbf$ from Section \ref{sec:motivation} as a (nonconvex) optimization problem over joint distributions. In Section \ref{sec:find} we introduce its convex relaxation. In Section \ref{sec:instances}, we solve problem $\Bbf$ by finding a solution to the KKT conditions of its relaxation that is a solution of the original nonconvex problem. In Section \ref{sec:witsen} we address the relaxation of the Witsenhausen problem.

\subsection{Problem formulation} \label{sec:formulation}
For the rest of the paper, we will consider the notation of Fig \ref{fig:comm}. Presently, we assume for simplicity that $\Sscr,\Xscr,\Yscr,\Sscrhat$ 
are finite. There are four random variables in the problem: $S,X,Y,\Shat$. Let $Q$ be their joint distribution and let `$Q_{\bullet}$'  denote the marginal of `$\bullet$'. Let lower case letters $s,x,y,\shat$ denote specific values of the respective upper case letters and `$Q_{V}(v)$' be abbrieviated as `$Q(v)$'. To be a distribution of $S,X,Y,\Shat$, $Q$ must satisfy the following `constraints'.
\begin{enumerate}
\item $Q$ must be a distribution.
 \item $Q$ must be consistent with the \textit{given marginals}. \ie 
\begin{equation*}
 Q_S(s)\equiv P_S(s), \quad \aur \quad Q_{Y|X}(y|x) \equiv P_{Y|X}(y|x).
\end{equation*} 
\item $Q$ must respect the \textit{information structure}, \ie $Q$ must satisfy Markovianity:
$$S\rightarrow X\rightarrow Y \rightarrow \Shat$$
\end{enumerate}
This notation says that $\Shat$ is independent of $X$ and $S$ given $Y$, and $Y$ is independent of $S$ given $X$. 
Therefore, a `feasible' distribution $Q$  is one that can be expressed as
\begin{equation}
Q(s,x,y,\shat)\equiv P_S(s)Q(x|s)P_{Y|X}(y|x)Q(\shat|y).  \label{eq:qdef}
\end{equation}  
We use $\Qscr$ to denote the set of all such $Q$.  
\begin{align}
\Qscr := \{Q |\  Q \in \Pscr(\Zscr), {\rm for\ which}\ \exists \ \non
Q_{X|S} \in \Pscr(\Xscr|\Sscr), \\Q_{\Shat|Y} \in \Pscr(\Sscrhat|\Yscr), \sthat Q \ {\rm satisfies} \
\eqref{eq:qdef} \},  \label{eq:qscrdef}
\end{align} 
where $\Pscr(\cdot)$ is the set of probability distributions on `$\cdot$' and we denote $\Zscr:= \Sscr \times \Xscr \times \Yscr \times \Sscrhat$. 

Any $Q \in \Qscr$ is parametrized by the kernels $Q_{X|S}$ and $Q_{\Shat|Y}$ (this parametrization has been studied in \cite{yuksel2012optimization}). 
In the context of Fig \ref{fig:comm}, these kernels can be thought of as ``random'' (single-letter) encoder and decoder, respectively. They together constitute 
what we call a \textit{random code}.
\begin{definition}
A random encoder is a conditional distribution $Q_{X|S}$ on $\Xscr$ given a symbol in $\Sscr$. A random decoder is a conditional distribution $Q_{\Shat|Y}$ on $\Sscrhat$ given a symbol in $\Yscr$. A pair of a random encoder and decoder is a \textit{random  code}.
An encoder (resp., decoder) is \textit{deterministic} if and only if there is a function $f$ (resp., $g$) such that 
$Q_{X|S}(x|s) = \I{f(s)=x}$ for all $x,s$ (resp., $Q_{\Shat|Y}(\shat|y)    =\I{g(y)=\shat}$ for all $y,\shat$). 
A pair of deterministic encoder and decoder is  a \textit{deterministic code}.
\end{definition}

The problem  that minimizes cost $\Ebb[ \kappa(S,X,Y,\Shat)]$ over random codes is our optimization-based formulation for $\Nbf$:
\begin{equation}
\ds \minimize{Q \in \Qscr} \left\langle\kappa,Q\right\rangle,   \tag{$\Nbf$}
\end{equation}
where $\langle\cdot,\cdot\rangle$ denotes the operation $\left\langle\kappa,Q\right\rangle := \sum_{z} \kappa(z)Q(z)$ 
and $z$ is a short-hand for a tuple $(s,x,y,\shat) \in \Zscr$. There is no abuse in denoting this problem as $\Nbf$ since it has the same optimal value as $\Nbf$ from Section \ref{sec:motivation}.
\begin{lemma}
The relaxed stochastic control problem over random codes (\ie $\Nbf$) has a solution that is a deterministic code. 
\end{lemma}
\begin{IEEEproof}
Write $\Nbf$ as
\begin{align*}
\min &\sum_{s,x,y,\shat}\kappa(s,x,y,\shat)P_S(s)Q_{X|S}(x|s)P_{Y|X}(y|x)Q_{\Shat|Y}(\shat|y) \\
&\subject\ {Q_{X|S} \in \Pscr(\Xscr|\Sscr), Q_{\Shat|Y} \in \Pscr(\Sscrhat|\Yscr) } 
\end{align*}
$\Nbf$ is thus a separably constrained bilinear optimization problem. It is a well known property of such a problem 
(see Exercise 4.25, \cite{bazaraa06nonlinear}) 
that it has a solution that is an extreme point; in this case a deterministic code.
\end{IEEEproof}
Without loss of generality we consider this $\Nbf$ as our stochastic control problem. \newline
\begin{remarkc}[Bilinearity:]
Note that the bilinearity of the objective of $\Nbf$ is a consequence of its information structure; it is not clear if this property will hold for joint distributions representing more exotic information structures.
\end{remarkc}

\def\Qscrdpi{\Qscr_{\rm DPI}}
\subsection{Convex relaxation} \label{sec:find}
We now formulate the convex relaxation of $\Nbf$. To do so, we first note the nonconvexity of this problem.
\begin{lemma} \label{lem:qscr}
$\Qscr$ is not a convex set.
\end{lemma}
\begin{IEEEproof}
Suppose $\Qscr$ is convex and consider distributions $Q^1,Q^2 \in \Qscr$ such that $Q^1_{X|S} \neq Q^2_{X|S}$ and $Q^1_{\Shat|Y} \neq Q^2_{\Shat|Y}$. Such $Q^1,Q^2$ exist so long as $\Xscr$ and $\Sscrhat$ are not both singletons.  Let $t \in (0,1)$. Then $q:=tQ^1 + (1-t)Q^2$ belongs to $\Qscr$ and satisfies \eqref{eq:qdef}. 
Consequently, we must have $q(\shat|y) \equiv q(\shat|s,x,y)$. Therefore using \eqref{eq:qdef},
\begin{align*}
 q(\shat|y)                        &= Q^1(\shat|y) + \frac{(1-t)Q^2(x|s)\left[ Q^2(\shat|y) -Q^1(\shat|y)\right] }{tQ^1(x|s) +
(1-t)Q^2(x|s)},
\end{align*}
for $s,x,y$ such that $q(s,x,y) >0$. The right hand side is independent of $x,s$  for all $\shat,y$ if and only if one of the following is true:
a) $Q^1(\shat|y) \equiv Q^2(\shat|y)$ or b) $\frac{(1-t)Q^2(x|s)}{tQ^1(x|s) +
(1-t)Q^2(x|s)} =c$, a constant, for all $x,s$. 
The first possibility is ruled out by the choice of $Q^1$ and $Q^2$. The second possibility implies 
$t(Q^1(x|s) - Q^2(x|s))\equiv 0$, which is ruled out since $t>0$. Thus $\Qscr$ is not convex.
\end{IEEEproof}

\begin{remarkc}
In the classical information structure $ Q $ would admit the decomposition
\[ Q(s,x,y,\shat)\equiv P_S(s)Q(x|s)P_{Y|X}(y|x)Q(\shat|s,x,y), \]
instead of \eqref{eq:qdef}. It is easy to see that the set of such $ Q $ is convex.
\end{remarkc}

Thus, $\Nbf$ is a nonconvex optimization problem. In general, one can only guarantee local 
optimality for such problems. A possible approach for finding a global solution is to solve its convex relaxation 
with the hope of finding a solution of the relaxation that is feasible for $\Nbf$. 

We consider the following relaxation: we minimize $\langle \kappa,Q\rangle$ over all $Q$ that satisfy the DPI and are consistent with the given source and channel distributions.  This problem is denoted $\Cbf_\Nbf$: 
\begin{equation}
 \minimize{Q \in \Qscrdpi} \left\langle \kappa, Q \right\rangle, \tag{$\Cbf_\Nbf$}
\end{equation} 
where $\Qscrdpi$ is the set:
\begin{align}
 \Qscr_{\rm DPI} = \{Q \in \Pscr(\Zscr) | Q_{S} = P_S, Q_{Y|X}=P_{Y|X} \\ \non
\aur \  I(Q_{X,Y}) \geq I(Q_{S,\Shat}) \}, 
\end{align} 
and $I(\cdot)$ is the mutual information under distribution `$\cdot$' and the log is the natural logarithm. Clearly, $\Qscr \subset \Qscr_{\rm DPI}$ (since Markovianity implies the DPI \eqref{eq:dpi}, but is not equivalent to it) and $\Qscrdpi$ is convex (this follows easily from the proof of Lemma \ref{lem:convexg} below).

To solve $\Cbf_\Nbf$, we write it as a mathematical progam using  additional 
variables $\{a(\shat|s) : \shat \in \Sscrhat, s \in \Sscr\}$,$\{b(x) : x \in \Xscr\}$.
 $$ \problemsmall{$\Cbf_{\Nbf}$}
         {Q,a,b}
         {\sum_{s,x,y,\shat}\kappa(s,x,y,\shat)Q(s,x,y,\shat)}
                                  {\small\hspace{-1.3cm}\begin{array}{r@{\ }c@{\ }l@{\quad}l@{\hspace{.55cm}}c@{}l@{}}
 \sum_{x,y}Q(z) &=& a(\shat|s) P_S(s), &\forall s,\shat, &:& \lambda^a(s,\shat), \\
 \sum_{s,y,\shat}Q(z) &=& b(x),  &\forall x, &:&\lambda^b(x),\\
 I(aP_S) & \leq  & I(P_{Y|X}b) & &:& \lambda,\\
 \sum_{s,\shat}Q(z) &=&P_{Y|X}(y|x)\sum_{s,y,\shat} Q(z) &\forall x,y, &:&\lambda^p(x,y),\\
 \sum_{\shat}a(\shat|s)    & = & \multicolumn{4}{l}{\hspace{-2mm}1,\ \forall s,\, :\hspace{-1mm}\mu^a(s)P_S(s), \quad\ \ \sum_x b(x) \hspace{-.5mm}=1,:\hspace{-0.5mm}\mu^b,}\\
 a(\shat|s) &\geq &\multicolumn{4}{l}{\hspace{-2mm}0,\ \forall s,\shat,\ :\hspace{-1mm} \nu^a(\shat|s)P_S(s),\ b(x) \geq\hspace{-1mm}0, \forall x, :\hspace{-1mm}\nu^b(x),}\\
 Q(z) &\geq & 0, &\forall z, &:&\nu(z).
 \end{array}} $$
Here $I(aP_S)$ is the mutual information of $S,\Shat$ under distribution $a(\shat|s)P_S(s)$ (and likewise $I(P_{Y|X}b)$). $\lambda^a,\lambda^b,\lambda,\mu^a,\mu^b,\nu^a,\nu^b$ and $\nu$ are (scaled) Lagrange multipliers corresponding to the specific constraints. Assume the existence of an interior (or Slater) point~\cite{luenberger97optimization}. Then $(Q,a,b)$ solves $\Cbf_\Nbf$ if and only if there exist Lagrange multipliers that 
satisfy the following system of KKT conditions 
\begin{align}
\tag{KKT$_{\Cbf_{\Nbf}}$} \lambda^a(s,\shat) &=  \lambda \log \frac{a(\shat|s)}{a(\shat)} + \mu^a(s) - \nu^a(\shat|s), \label{eq:kktcbfnbf} \\ 
 \lambda^b(x) &=- \lambda  D\left( P_{Y|X}(\cdot|x)||b_Y(\cdot) \right) + \lambda  +\mu^b-\nu^b(x),  \non \\
\nu(z) &= \kappa(z) + \lambda^a(s,\shat) + \lambda^b(x) + \lambda^p(x,y)\non \\ 
& \textstyle \qquad -\sum_{\bar{y}} \lambda^p(x,\bar{y})P_{Y|X}(\bar{y}|x),  \tag{$*$}\label{eq:nu} \\
0\leq \lambda      & \perp I(P_{Y|X}b)-I(aP_S) \geq 0, \non \\
\left\langle\nu,Q\right\rangle &=0, \ \nu(\cdot) \geq 0, Q(\cdot) \geq 0 \textstyle\sum_{\shat}a(\shat|s)    =  1,\ \ \forall z \in \Zscr, \non  
\end{align}
where we have used, $\frac{\del I(aP_S)}{\del a(\shat|s)} = P_S(s) \log \frac{a(\shat|s)}{a(\shat)},$
$\frac{\del I(P_{Y|X}b)}{\del b(x)} =  D\left( P_{Y|X}(\cdot|x)||b_Y(\cdot) \right)-1$ and $b_Y(\cdot)$ denotes $\sum_x P_{Y|X}(y|x)b(x)$.

To solve $\Nbf$ for optimality or inverse optimality by the logic of the convex relaxation, we have to tackle \eqref{eq:kktcbfnbf}. The hardness of this depends on structure of $\kappa$ and the choice of the convexification, as we explain below.\newline
\begin{remarkc}[DPI and the structure of $\kappa$:] \label{rem:kappastr} Perhaps the most difficult part of \eqref{eq:kktcbfnbf} is \eqref{eq:nu}, which contains all variables involved and the left hand side $\nu$  has to be both nonnegative and complementary to $Q$. One way of dealing with this would be to set $\nu(\cdot) \equiv 0$. Then \eqref{eq:nu} necessitates that $\kappa$ be of the form
\begin{equation}
k(s,x,y,\shat) \equiv f_1(x) + f_2(s,\shat) + f_3(x,y), \label{eq:kappastr}
\end{equation}
where the functions $f_1,f_2,f_3$ are determined by the $Q$ that solves $\Cbf_\Nbf$. The communication setting, where $\kappa(z) \equiv \delta(s,\shat) + \rho(x)$, agrees naturally with this structure and therefore taking $\nu(\cdot)\equiv0$ may solve \eqref{eq:kktcbfnbf}. Furthermore, \textit{inverse} optimality is easy for this structure (in fact, $\nu(\cdot)\equiv 0$ gives the same form of inverse optimal $\delta$ and $\rho$ as \eqref{eq:dchar}). Problem $\Bbf$, on the other hand, has $\kappa(z) \equiv \delta(s,\shat) + \rho(x) + \tau(x,s)$ and does not agree directly with \eqref{eq:kappastr}. In our solution of $\Bbf$, we will not take $\nu(\cdot)\equiv0$. The Witsenhausen problem has  $\kappa(z) \equiv \zeta(x,\shat) + \tau(x,s)$ and it is even harder to reconcile with \eqref{eq:kappastr} (see Theorem \ref{thm:witsen}). The mismatch between $\kappa$ and \eqref{eq:kappastr} makes inverse optimality for these problems hard too.

The structure of the right hand side of \eqref{eq:kappastr} is due the DPI constraint. Other convexifications would yield other structures and may allow greater ease in solving \eqref{eq:kktcbfnbf}. What we have developed so far is thus a \textit{general framework} of which the DPI-based convexification is a particular application.\end{remarkc}

In the following sections, we show that this view of the DPI as a convexification tool is sufficient for solving  $\Bbf$.  In Section \ref{sec:gastparsol} we solve a variant of $\Bbf$ with $\kappa=\delta+\rho$. In Section \ref{sec:bansalsol} we solve $\Bbf$. As argued in the above remark, when $ \kappa = \delta + \rho, $ we take $ \nu \equiv 0, $ whereas for $ \Bbf, $ our solution uses $ \nu \not\equiv 0. $ 
\subsection{Solution of specific instances} \label{sec:instances}

\subsubsection{Solution with Gastpar-type cost functions}\label{sec:gastparsol}
In this section consider a special case of $\Nbf$, denoted $\Gbf$, where $\kappa(z) \equiv \delta(s,\shat) + \rho(x)$:  \vspace{-30pt}
\begin{multicols}{2}
\begin{equation}
\minimize{Q \in \Qscr} \left\langle \delta + \rho,Q \right\rangle, \tag{$\Gbf$}
\end{equation}\break
\begin{equation}
\minimize{Q \in \Qscrdpi} \left\langle \delta + \rho,Q\right\rangle. \tag{$\Cbf_\Gbf$}
\end{equation}
\end{multicols}\vspace{-15pt}
\noindent and $\Cbf_{\Gbf}$ is the convex relaxation of $\Gbf$. 
We first write $\Cbf_\Gbf$ as a mathematical program in a simplier form than $\Cbf_\Nbf$. Notice that the variable `$Q$' in $\Cbf_\Nbf$ can be dropped since the objective of $\Cbf_\Gbf$ can be expressed in terms of $a,b$ alone, 
and for any $(a,b)$ in the set $\{(a,b) | a \in \Pscr(\Sscrhat|\Sscr), b \in \Pscr(\Xscr), I(aP_S) \leq I(P_{Y|X}b) \},$ there exists $Q$ satisfying the constraints in $\Cbf_\Nbf$ (take $Q(z) \equiv a(\shat|s)P_S(s)b(x)P_{Y|X}(y|x)$). 
An equivalent form of $\Cbf_\Gbf$ is thus the following.
$$ \problemsmall{$\Cbf_{\Gbf}$}
        {a,b}
        {\sum_{s,\shat}\delta(s,\shat)a(\shat|s)P_S(s) + \sum_{x}\rho(x)b(x)}
                                 {\hspace{-1.1cm}\begin{array}{r@{\, }c@{\, }l@{}c@{}}
\sum_{\shat}a(\shat|s)    & = & 1,  \forall s,\quad :\hspace{-1mm}\mu^a(s)P_S(s), & \sum_{x}b(x)             =  1, :\hspace{-1mm}\mu^b,\\ 
I(aP_S)                  & \leq  & I(P_{Y|X}b), \qquad :\lambda\\  
 a(\shat|s)  \geq &0,& \multicolumn{2}{l}{\hspace{-2mm}  \forall s,\shat\ \ :\hspace{-1mm}\nu^a(\shat|s)P_S(s),  \ \  b(x)     \hspace{-1mm}\geq  0,\, \forall x,\hspace{-.5mm} :\hspace{-1mm}\nu^b(x)}
        \end{array}} $$
\begin{lemma} \label{lem:convexg}
$\Cbf_\Gbf$ is a convex optimization problem.
\end{lemma}
\begin{IEEEproof}
The objective is linear. All constraints except the third are linear. It is well known~\cite{cover06elements} that for fixed $P_S$,  $I(aP_S)$ is convex in $a$ and for fixed $P_{Y|X}$, $I(P_{Y|X}b)$ is concave in $b$, whereby the feasible region of is convex. 
\end{IEEEproof}

Let $\mu^a,\mu^b, \lambda,\nu^a,\nu^b$ be the Lagrange multipliers of the constraints of problem $\Cbf_\Gbf$.  
$a,b$ are optimal for $\Cbf_\Gbf$ if and only if together with the Lagrange multipliers, they satisfy,
\begin{align}
 \delta(s,\shat)  =&-  \lambda  \log \frac{a(\shat|s)}{a(\shat)} - \mu^a(s) + \nu^a(\shat|s) \tag{KKT$_{\Cbf_\Gbf}$}\label{eq:kktnbf}\\ 
 \rho(x)             =& \lambda  D\left( P_{Y|X}(\cdot|x)||b_Y(\cdot) \right) -\lambda- \mu^b + \nu^b(x) \non \\ 
0\leq \lambda       \perp& I(P_{Y|X}b)-I(aP_S) \geq 0 \non\\
\textstyle \sum_{\shat}a(\shat|s) &= 1,  \left\langle a(\cdot|s), \nu^a(\cdot|s)\right\rangle    = 0, \nu^a(\shat|s)\geq 0, \non  \\
\textstyle \sum_{x}b(x)             =&  1, \left\langle b,\nu^b\right\rangle  =0, \nu^b(\cdot) \geq 0, b(\cdot) \geq 0,  a(\shat |s)\geq 0,\non
\end{align}
$\forall z\in \Zscr$. By making this simplification, we have effectively taken $\nu(\cdot) \equiv 0$, $\lambda^p(\cdot)\equiv0$ in \eqref{eq:kktcbfnbf} and in \eqref{eq:nu} we are separately setting terms in $(s,\shat)$ and $x$ to be zero. 

We are now prepared to demonstrate our goal which was to find a solution of $\Cbf_\Gbf$ that is feasible for $\Gbf$, in order to arrive at a solution to $\Gbf$. Although we have taken $\Zscr$ to be finite, we assume that \eqref{eq:kktnbf} holds for this problem too. The arguments used are essentially geometric and they can be generalized with appropriate technical assumptions. We abbreviate $\xi_0(\gamma_0) :=\frac{\sigma_0^2\sigma_w^2}{({\gamma_0}^2\sigma_0^2 + \sigma_w^2)^2}$, $\xi_1(\gamma_0):=\frac{\gamma_0\sigma_0^2}{{\gamma_0}^2\sigma_0^2 + \sigma_w^2}$ and $\bar{P}(\gamma_0) := \frac{\gamma_0^2\sigma_0^2+ \sigma_w^2}{\gamma_0^2\sigma_0^2}$. 
\begin{theorem} \label{thm:gastpar}
Consider problem $\Bbf$  with $s_{01}=0$. \ie consider $\Gbf$ where $S \sim \Nscr(0,\sigma_0^2)$, $Y = X + w,$ $w \sim \Nscr(0,\sigma_w^2)$ and independent of $S$, $\delta(s,\shat) \equiv (\shat-s)^2$ and $\rho(x) \equiv k_0x^2$. Assume that a solution of $\Cbf_\Gbf$ is characterized by the KKT conditions \eqref{eq:kktnbf}. 
Then the PDF given by $(S,\gamma_0^* S,Y,\gamma_1^*Y)$ where $\gamma_0^*$ satisfies $k_0 = \xi_0(\gamma_0^*)$ and $\gamma_1^*=\xi_1(\gamma_0^*)$ solves the convex relaxation $\Cbf_\Gbf$. Furthermore, it is a solution of $\Gbf$.
\end{theorem}
\begin{IEEEproof} We parametrize $X = \gamma_0S$ and $\Shat = \gamma_1 Y$ and find 
$\gamma_0$ and $\gamma_1$ to solve \eqref{eq:kktnbf}. The distribution obtained would be feasible for $\Gbf$ since $S \rightarrow \gamma_0 S \rightarrow Y \rightarrow \gamma_1 Y$, and would thus be a solution of $\Gbf$. With these constants we get (details can be found in \cite{gastpar03thesis}, p.\ 51)
\begin{align}
 D(P_{Y|X}(\cdot|x) || b_Y(\cdot)) &= c_1+ \frac{x^2}{2\gamma_0^2\sigma_0^2\bar{P}} \label{eq:gaussd}\\
 \log\frac{a(\shat|s)}{a(\shat)} = \frac{-1}{2\gamma_1^2\sigma_w^2\bar{P}}&\left[\shat - \bar{P}\gamma_0\gamma_1s\right]^2+ c_2(s) \label{eq:gaussrho}
 \end{align}
where $\bar{P} = \bar{P}(\gamma_0)$, $c_1$ independent of $x$ 
and $c_2$ is independent of $\shat$. 
To solve system \eqref{eq:kktnbf}, take  $\nu^b(x),\nu^a(\shat|s) \equiv 0$, $\mu^b = \lambda c_1-\lambda$ and $\mu^a(s) \equiv -\lambda c_2(s)$. The first two equations in \eqref{eq:kktnbf} reduce to the identities 
\begin{align*}
k_0x^2 \equiv \lambda \frac{x^2}{2\gamma_0^2\sigma_0^2\bar{P}}, \quad (\shat-s)^2 \equiv \frac{\lambda}{2\gamma_1^2\sigma_w^2\bar{P}}\left[\shat - \bar{P}\gamma_0\gamma_1s\right]^2.
\end{align*}
Comparing coefficients, we get $k_0 = \xi_0(\gamma_0), \gamma_0\gamma_1\bar{P}=1$ and  $\lambda=2\gamma_1^2\sigma_w^2\bar{P}.$ Simplifying, we get $\gamma_0=\gamma_0^*$, where $\gamma_0^*$ satisfies $k_0 = \xi_0(\gamma_0^*)$ and $\gamma_1 = \gamma_1^* = \xi_1(\gamma_0^*)$. To show that this solves \eqref{eq:kktnbf}, we need to show complementarity slackness of the DPI constraint. The values $\gamma_0^*$ and $\gamma_1^*$ imply $\lambda>0$. Therefore we need to show that equality holds in the DPI. This follows by noting that $\Ebb[X^2]= {\gamma_0^*}^2\sigma_0^2$ and $\Ebb[(\Shat - S)^2]  = \frac{\sigma_w^2\sigma_0^2}{{\gamma_0^*}^2\sigma_0^2 + \sigma_w^2}$, which means that the inequalities \eqref{eq:bbbound1}-\eqref{eq:bbbound2} are tight with $P$ taken as ${\gamma_0^*}^2\sigma_0^2$. Thus, the distributions induced by taking $X =\gamma_0^*S$ and $\Shat = \gamma_1^*Y$ solve \eqref{eq:kktnbf}  and are  optimal for $\Cbf_\Gbf$. Since they satisfy Markovianity, they also solve $\Gbf$.
\end{IEEEproof}

\begin{remarkc}[Gaussian test channel:] That these calculations look familiar and that the gains $\gamma_0^*$ and $\gamma_1^*$ are the well-known codes for the Gaussian test channel is not surprising; when the objective is a function of $a,b$ alone, \eqref{eq:kktcbfnbf} reduces to \eqref{eq:kktnbf}, the right hand side of which has the same form as Gastpar's \eqref{eq:dchar}. But recall that our calculations do not pertain to a communication problem. They are for the solution of a convex relaxation of a nonconvex \textit{optimization} problem which was convexified using the DPI. 
\end{remarkc}

$\gamma_0^*$ and $\gamma_1^*$ agree with the solution of Bansal and Ba\c{s}ar for $s_{01}=0$. But the important consequence, which is part of the message of this paper, is that, problem $\Bbf$ (with $s_{01}=0$) shares its solution with a specific convex relaxation. We emphasize this through Theorem \ref{thm:gastparemph} below.
\begin{theorem} \label{thm:gastparemph}
Let the setting and the assumptions of Theorem \ref{thm:gastpar} hold.  
The optimal value of problem $\Bbf$ with $s_{01}=0$ equals the optimal value of its 
convex relaxation with $\Qscrdpi$. Furthermore, any solution of $\Bbf$ with $s_{01}=0$ solves its convex relaxation with $\Qscrdpi$.
\end{theorem}

\subsubsection{Solution of the Bansal-Ba\c{s}ar problem}\label{sec:bansalsol}
We now return to problem $ \Bbf $. 
$\Bbf$ has the following structure: $\kappa(z) \equiv \delta(s,\shat) + \rho(x) + \tau(x,s),$ where $\tau(x,s) \equiv k\sqrt{\rho(x)}\tau'(s)$.  We provide two proofs for the solution of $ \Bbf. $ The first proof uses a reasoning similar to that used by Bansal and Ba\c{s}ar --  specifically, the use of the Cauchy-Swartz inequality for bounding the term $ \tau(x,s).$ In our second proof (which is our new proof), we do not use this inequality and argue directly. 

Let $\Cbf_\Bbf$ be the convex relaxation of problem $ \Bbf $. 

\noindent\begin{tabular}{p{.23\textwidth}p{.23\textwidth}}
{\begin{align}
\minimize{Q \in \Qscr} \left\langle \delta + \rho + \tau, Q \right\rangle \tag{$\Bbf$}
\end{align}}
&
{\begin{align}
\hspace{-.5cm}\minimize{Q \in \Qscrdpi} \left\langle \delta + \rho +\tau, Q \right\rangle \tag{$\Cbf_\Bbf$}
\end{align}}
\end{tabular}
By Cauchy-Schwartz inequality, we get $k\sum_{x,s}\sqrt{\rho(x)}\tau'(s)P_S(s)Q(x|s) \geq -\alpha (\sum_x\rho(x)b(x))^{1/2},$ where $\alpha = |k|\sqrt{\sum_s (\tau'(s))^2P_S(s)}.$ The use of the Cauchy-Schwartz inequality allows for the construction of an auxilliary problem, denoted $\Bbf'$, where $\langle\tau,Q\rangle$ is replaced with $-\alpha\sqrt{\langle \rho,Q\rangle}$
\begin{equation}
\minimize{Q \in \Qscr} \left\langle \delta + \rho, Q \right\rangle  - \alpha  \sqrt{\left\langle \rho, Q\right\rangle}. \tag{$\Bbf'$}
\end{equation} 
The convex relaxation of $\Bbf'$ using $\Qscrdpi$ is the following,
\begin{equation}
 \minimize{Q \in \Qscrdpi} \left\langle \delta + \rho, Q \right\rangle  - \alpha  \sqrt{\left\langle \rho, Q\right\rangle}. \tag{$\Cbf_{\Bbf'}$}
\end{equation} 
We will use $\Cbf_{\Bbf'}$ to solve $\Bbf$. Note that $\Cbf_{\Bbf'}$ is \textit{not} a special case of $\Cbf_\Nbf$ because its objective is not linear in $Q$. But it can  be written as a mathematical program, as $\Cbf_\Gbf$ was, in terms of the variables $a,b$. The resulting formulation of $\Cbf_{\Bbf'}$ has the same constraints as $\Cbf_\Gbf$ and its objective is 
$$\sum_{s,\shat}\delta(s,\shat)P_S(s)a(\shat|s) + \sum_{x}\rho(x)b(x) -\alpha (\sum_x \rho(x)b(x))^{1/2}.$$
\begin{lemma}
$\Cbf_{\Bbf'}$ is a convex optimization problem.
\end{lemma}
\begin{IEEEproof}
It suffices to show $\psi(b) :=-\sqrt{\sum_x \rho(x)b(x)}$ is a convex function of $b$ (cf. Lemma \ref{lem:convexg}). 
The Hessian of  $\psi$, $\nabla^2\psi(b) = \frac{1}{4}\frac{1}{(\sum_x \rho(x)b(x))^{3/2}} \vec{\rho} \vec{\rho}^\top$,  where $\vec{\rho} :=(\rho(x))_{x \in \Xscr}$ is the vector of values of $\rho$. Clearly, $\nabla^2 \psi(b) \succeq 0$ for all $b$. 
\end{IEEEproof}
Consequently, the KKT conditions of $\Cbf_{\Bbf'}$ characterize its solution. We will use this to solve $\Bbf$.
\begin{theorem} \label{thm:bb}
Consider the problem of Bansal and Ba\c{s}ar. \ie, consider problem $\Bbf$ where $S \sim \Nscr(0,\sigma_0^2)$, $Y = X + w$ where $w \sim \Nscr(0,\sigma_w^2)$ and independent of $S$, $\delta(s,\shat) \equiv (\shat-s)^2$, $\rho(x) \equiv k_0x^2$ and $\tau(x,s) \equiv s_{01}xs$. Suppose the solution of $\Cbf_{\Bbf'}$ is characterized by its KKT conditions. Then the PDF  $(S,\gamma_0^{*} S,Y,\gamma_1^{*}Y)$ where $\gamma_1^{*} = \xi_1(\gamma_0^{*})$
and $\gamma_0^{*}$ is a positive solution of 
\begin{equation}
 (2k_0\gamma_0^*\sigma_0 - |s_{01}|\sigma_0)({\gamma_0^*}^2 \sigma_0^2 + \sigma_w^2)^2 =    2\gamma_0^*\sigma_0^3\sigma_w^2, \label{eq:bbgamma}
\end{equation} 
solves $\Bbf'$ and its convex relaxation $\Cbf_{\Bbf'}$. Furthermore, the PDF $(S,\gamma_0^{**} S,Y,\gamma_1^{**}Y)$ where $\gamma_1^{**} = \xi_1(\gamma_0^{**})$ and $\gamma_0^{**}=-\sgn(s_{01})\gamma_0^*$ solves $\Bbf$.
\end{theorem}
\begin{IEEEproof}
By Cauchy-Schwartz inequality and since $\Qscr \subset \Qscrdpi$ we have 
$$\OPT(\Bbf) \geq \OPT(\Bbf') \geq \OPT(\Cbf_{\Bbf'}), $$
where $\OPT(\cdot)$ is the optimal value of `$\cdot$'. 
We first show that the second inequality is tight 
by finding a solution of $\Cbf_{\Bbf'}$ that is feasible for $\Bbf'$. Then we will construct a feasible point of 
$\Bbf$ for which the objective of $\Bbf$ equals $\OPT(\Bbf')$. 

The KKT conditions of $\Cbf_{\Bbf'}$ are similar to those of $\Cbf_\Gbf$, except for the equation
\begin{align}
\rho(x)(1 - \alpha/(2R))            = \lambda  D\left( P_{Y|X}(\cdot|x)||b_Y(\cdot)
 \right) -\lambda - \mu^b + \nu^b(x), \non 
\end{align}
where $R:= \sqrt{\sum_t \rho(t)b(t)}$. Here, $\alpha = \frac{|s_{01}|\sigma_0}{\sqrt{k_0}}$. 
We postulate that $X = \gamma_0S, \Shat = \gamma_1Y$ to solve the KKT conditions of $\Cbf_{\Bbf'}$. Assume $\gamma_0>0$, so $R = \sqrt{k_0}\gamma_0\sigma_0$. As in Theorem \ref{thm:gastpar}, put $\nu^b(x),\nu^a(\shat|s) \equiv 0$, $\mu^b = \lambda c_1 - \lambda$, $\mu^a(s) \equiv -\lambda c_2(s)$ and  compare coefficients to get $k_0 = \frac{\lambda}{2(1-\frac{\alpha}{2R})(\gamma_0^2 \sigma_0^2 + \sigma_w^2)},$ $\gamma_0\gamma_1\bar{P}=1$ and $\lambda=2\gamma_1^2\sigma_w^2\bar{P}$. Simplifying gives $\gamma_0=\gamma_0^*$, where $\gamma_0^*$ satisfies \eqref{eq:bbgamma} and $\gamma_1 = \gamma_1^*=\xi_1(\gamma_0^{*})$. 
The equality in DPI follows as in Theorem \ref{thm:gastpar}. Therefore, the distribution $(S,\gamma_0^*S,Y,\gamma_1^*Y)$ is optimal for $\Cbf_{\Bbf'}$ and, by Markoviantiy, for $\Bbf'$. Thus we get that $$\OPT(\Bbf') = \frac{\sigma_w^2\sigma_0^2}{{\gamma_0^*}^2\sigma_0^2 + \sigma_w^2} + k_0{\gamma_0^*}^2\sigma_0^2 -|s_{01}|\gamma_0^*\sigma_0^2.$$ 

Now consider the PDF  $Q^{**}$ defined by $(S,\gamma_0^{**}S,Y,\gamma_1^{**}Y)$, where $\gamma_0^{**}\in \Real$ and $\gamma_1^{**} = \xi_1(\gamma_0^{**})$. It is easy to check if $\gamma_0^{**} = -\sgn(s_{01})\gamma_0^*$, the objective of $\Bbf$ evaluated at $Q^{**}$ 
equals $\OPT(\Bbf')$. Thus $Q^{**}$ solves $\Bbf$. 
\end{IEEEproof}

Note that \eqref{eq:bbgamma} agrees with \eqref{eq:bbpstar} of Bansal and Ba\c{s}ar with $P^* = \gamma_0^*\sigma_0$ and the $\gamma_0^{**}$ agrees with the optimal controller they obtain in \cite{bansal87stochastic}. Although the above result did not utilize the convex relaxation of $ \Bbf $, it is indeed true that a solution of $\Bbf$ also solves its convex relaxation with $\Qscrdpi$.
Below we emphasize this.
\begin{theorem} \label{thm:bbequal}
Let the setting and assumptions of Theorem \ref{thm:bb} hold. Any solution of $\Bbf$ solves its relaxation $\Cbf_\Bbf$ with $\Qscrdpi$. \ie the optimal value of $\Bbf$ equals the optimal value of $\Cbf_\Bbf$. In fact, $\OPT(\Bbf) = \OPT(\Bbf') = \OPT(\Cbf_\Bbf) = \OPT(\Cbf_{\Bbf'}).$
\end{theorem}
\begin{IEEEproof}
We only need to prove the last sequence of equalities. Since $\Qscr \subset \Qscrdpi$, $\OPT(\Bbf)\geq \OPT(\Cbf_\Bbf)$ and by  Cauchy-Schwartz inequality, $\OPT(\Cbf_\Bbf) \geq \OPT(\Cbf_{\Bbf'})$. But Theorem \ref{thm:bb} shows $\OPT(\Bbf) =\OPT(\Cbf_{\Bbf'})=\OPT(\Bbf')$. Combining this gives the result.
\end{IEEEproof}

We now show that as a consequence of Theorem \ref{thm:bbequal}, any solution of $\Bbf$ must satisfy equality in the DPI. For this we need the following lemma.
\begin{lemma} \label{lem:cgprime} Let the setting and assumptions of Theorem \ref{thm:bb} hold and  let $Q^*$ be a solution of $\Bbf$. Then $Q^*$ also solves 
\begin{equation}
\minimize{Q \in \Qscrdpi}\langle \delta + (1-\alpha/(2R^*))\rho,Q\rangle \tag{$\Cbf_{\Gbf'}$}
\end{equation}
where $R^* = \sqrt{\langle \rho,Q^*\rangle}$.
\end{lemma}
\begin{IEEEproof}
From Theorem  \ref{thm:bbequal}, $Q^*$ solves $\Cbf_{\Bbf'}$. The proof now follows easily by comparison of the 
KKT conditions of $\Cbf_{\Bbf'}$ and $\Cbf_{\Gbf'}$.
\end{IEEEproof}
\begin{theorem} \label{thm:dpiequal}
Let the setting and assumptions of Theorem \ref{thm:bb} hold and  let $Q^*$ be a solution of $\Bbf$. Define the rate-distortion function $R(\cdot)$ and capacity cost function $C(\cdot)$ as in \eqref{eq:rd}-\eqref{eq:cp} for distortion $\delta$ and cost $(1-\alpha/(2R^*))\rho$, where $R^* = \sqrt{\langle \rho,Q^*\rangle}$. 
Suppose for each $P \geq 0$ there exists $D \geq 0$ such that $C(P) = R(D)$.  Then the equality holds in the DPI under the distribution $Q^*$.
\end{theorem}
\begin{IEEEproof} By Lemma \ref{lem:cgprime}, $Q^*$ solves $\Cbf_{\Gbf'}$. 
Let $a=Q^*_{\Shat|S},b=Q^*_{X}$ and for sake of contradiction, suppose $I(aP_S) < I(P_{Y|X}b)$. Clearly, for $D = \Ebb_a[\delta]$ and $P=\Ebb_b[(1-\alpha/(2R^*))\rho]$, $ R(D) \leq I(aP_S) < I(P_{Y|X}b) \leq C(P)$. 
By the hypothesis, there exists $D'$ with $R(D')=C(P)>0$ and since $R(\cdot)$ is a strictly decreasing when it is positive, 
$D'<D$. Let $a'$ be a distribution on $\Shat|S$ such that $D'=\Ebb_{a'}[\delta]$. Then $Q'(z) \equiv P_S(s)a'(\shat|s)P_{Y|X}(y|x)b(x)$ is feasible for $\Cbf_{\Gbf'}$ and has a lower value than $Q^*$. A contradiction.
\end{IEEEproof}
\begin{remarkc}[DPI in problem $\Bbf$:]
It is an open problem to know if there a way of solving $\Bbf$ that does not use the DPI~\cite{basar08variations}. The above theorem says that, except in a degenerate case, equality in the DPI is necessary for optimality of $ \Bbf $. Thus any other solution method must imply this equality. Note that this is true even for solutions of $\Bbf$ that are not `linear' (where $X = \gamma_0 S$, $\Shat = \gamma_1 Y$)  or deterministic. 
\end{remarkc}

The original proof of Bansal and Ba\c{s}ar used Cauchy-Schwartz inequality (cf. Section \ref{sec:bb}). We used this  inequality in the construction of the auxilliary problem $\Bbf'$. We now present an alternative proof by directly solving the convex relaxation $\Cbf_\Bbf$. But unlike $\Cbf_{\Bbf'}$ or $ \Cbf_\Gbf $, the objective of $\Cbf_\Bbf$ cannot written in terms of $a,b$ alone and we have to deal with the harder system \eqref{eq:kktcbfnbf}.

\begin{IEEEproof}[Proof of Theorem \ref{thm:bb} (without Cauchy-Schwartz)]
We will solve \eqref{eq:kktcbfnbf} for $\Cbf_\Bbf$. Let $X \sim \gamma_0S$, $\Shat \sim \gamma_1Y$ and $Q$ be the 
resulting PDF. In \eqref{eq:kktcbfnbf}, use \eqref{eq:gaussd}-\eqref{eq:gaussrho}, and put $\gamma_1=\xi_1(\gamma_0)$, $\lambda = 2\gamma_1^2\sigma_w^2\bar{P}$, $\mu^b = \lambda c_1 - \lambda$, $\nu^a(\cdot),\nu^b(\cdot),\lambda^p(\cdot) \equiv0,$
and $\mu^a(s) \equiv -\lambda c_2(s) + \bar{\mu}^a(s)$, where $\bar{\mu}^a(s)$ will be determined later. This reduces \eqref{eq:nu} to $$\nu(z) \equiv k_0x^2 - \xi x^2 + s_{01}xs + \bar{\mu}^a(s),$$
where $\xi := \xi_0(\gamma_0).$  If $s_{01} =0$, setting $k_0 =\xi$ 
and $\bar{\mu}^a(\cdot)\equiv 0$ solves the problem. So we assume $s_{01} \neq 0$. We limit our choice of $\gamma_0$ to those which satisfy 
$k_0 >\xi_0(\gamma_0)$ and first find $\bar{\mu}^a(s)$ so as to satisfy $\nu(\cdot) \geq 0$. Completing the squares gives 
$$\nu(z) \equiv   \left(  \sqrt{k_0 - \xi}x + \frac{s_{01}}{2\sqrt{k_0 - \xi}}s\right)^2 - \frac{s_{01}^2s^2}{4(k_0-\xi)} + \bar{\mu}^a(s). $$
Put $\bar{\mu}^a(s) \equiv \frac{s_{01}^2s^2}{4(k_0-\xi)}$, so that now $\nu(\cdot) \geq 0$. The complementarity slackness of the DPI constraint follows as before. So we only need to satisfy the complementarity of $\nu(\cdot)$ and $Q(\cdot)$. Observe 
that $\nu(z) >0$ if and only if $x \neq -\frac{s_{01}s}{2(k_0-\xi)}$, whereas for $x \neq \gamma_0 s$, the 
PDF $Q(z)=0$. Thus, it suffices to 
take $\gamma_0 = -\frac{s_{01}}{2(k_0-\xi)}$ to satisfy the complementarity of $\nu(\cdot)$ and $Q(\cdot)$. 
Therefore taking $\gamma_0=\gamma_0^{**}=-\sgn(s_{01})\gamma_0^*$, where $\gamma_0^*$ is a positive solution of 
\eqref{eq:bbgamma}, satisfies $\gamma_0 = -\frac{s_{01}}{2(k_0-\xi)}$. The positivity of $\gamma_0^*$ ensures 
that $k_0 >\xi$ if $s_{01} \neq 0$. This completes the proof.
\end{IEEEproof} 
\begin{remarkc}[Separability of the cost:]
 Although the objective of $\Bbf$ cannot be written in $a,b$ alone, using Cauchy-Schwartz inequality, a problem $\Bbf'$ with $\OPT(\Bbf')=\OPT(\Bbf)$, and whose objective could be written in $a,b$ was constructed. For a problem with this structure, results from communication theory apply more directly (due to the structure of the DPI) and Bansal and Ba\c{s}ar were able to leverage them. The optimization approach, on the other hand, can address $\Bbf$ more directly and through a more general framework, without the need for Cauchy-Schwartz inequality. \end{remarkc}

In summary we have shown that $\Bbf$ can be solved by  a particular application of the broader concept of convex relaxation where the DPI was used for constructing the convex relaxation.

\section{INVERSE OPTIMAL CONTROL} \label{sec:ioc}
We now consider inverse optimal control and show how it may be addressed  from the standpoint of convex relaxation. The specific problem we are interested in is problem $\bf G$, where $\kappa(z) \equiv \delta(s,\shat) + \rho(x)$, although the overall framework applies to any case of $\Nbf$.  

The thrust of this section is in showing that convex relaxations allow one to obtain a wider range of inverse optimal cost functions for $\Gbf$ than are obtainable through the communication theoretic argument. To explain this, 
consider the convex relaxation $\Cbf_\Gbf$ of $\Gbf$ and recall the system \eqref{eq:kktnbf}. 
We first clarify the distinction between inverse optimality for $\Gbf$ and that for Gastpar's problem which is  characterized by \eqref{eq:dchar}.

\begin{remarkc}[Inverse optimal control:] \label{rem:ioc}
Suppose we are interested in inverse optimality of $\delta$ and $\rho$ for problem $\Gbf$. Nonconvexity of $\Qscr$ makes this hard, but the normal cone of $\Qscrdpi$ provides a subset of the inverse optimal functions (cf. Section \ref{sec:message}). The right-hand side of the system \eqref{eq:kktnbf} is precisely this subset. Although the form of $\delta$ and $\rho$ in \eqref{eq:kktnbf} are the same as Gastpar's in \eqref{eq:dchar} (except for minor differences), ours are obtained for an altogether  different notion of optimality. Gastpar has inverse optimality \textit{over deterministic codes of arbitrary block length}; ours are for  inverse optimality over $\Qscr$, \ie, over single-letter random codes. For inverse optimality for $\Gbf$, \eqref{eq:kktnbf} is  \textit{not} a necessary condition (unlike for optimality over arbitrary block lengths where \eqref{eq:dchar} is also necessary). They have coincided because \textit{mutual information serves the dual role} of characterizing Gastpar's optimality and of convexifying $\Qscr$. The results of the previous sections demonstrate that for solving $ \Bbf $ or $ \Gbf, $ only the latter role is important.

 There are minor differences between \eqref{eq:kktnbf} and \eqref{eq:dchar}. 
 In \eqref{eq:dchar} the constants $c_1,c_2$ are allowed to be distinct, but in \eqref{eq:kktnbf} they are equal ($ =\lambda $). This is because Gastpar considers \textit{Pareto} optimality; a Pareto optimal code minimizes $\delta + \alpha \rho$ for some $\alpha>0$. We have fixed $\alpha=1$. Gastpar requires $c_1,c_2>0$, whereas our $\lambda$ may be zero. Gastpar has an additional condition that  $R(\Delta^*)$ and $C(\Gamma^*)$ cannot be held fixed while reducing the optimal $\Delta^*$ and $\Gamma^*$. This is equivalent, in our problem, to the sensitivity of the DPI constraint, which implies the strict positivity of $\lambda$. Equality in the DPI (analogous to $R(\Delta^*)=C(\Gamma^*)$) is a necessary condition for optimality in $\Cbf_\Gbf$ too (this follows from the proof of Theorem \ref{thm:dpiequal}). Finally, Gastpar does not have a term corresponding to $\nu^a(\shat|s)$; he has deliberately dropped it~\cite{gastpar03thesis} since when $\nu^a(\shat|s)>0$, $\log \frac{a(\shat|s)}{a(\shat)}$ is undefined. 
\end{remarkc}

We have chosen to focus on problem $\Gbf$ because for this problem the inverse optimal cost functions can be read off directly from the KKT conditions. In general for problem $\Nbf$, one would have to obtain the cost functions by solving the system \eqref{eq:kktcbfnbf}. The difficulty of doing this was articulated in Remark \ref{rem:kappastr}.

The main implication of Remark \ref{rem:ioc} is as follows. The convexifcation using $\Qscrdpi$ is only one convexifcation of $\Qscr$ that provided a specific subset of the inverse optimal cost functions of $\Gbf$, namely the one that coincided with  Gastpar's inverse optimal cost functions. One could potentially employ other convexifications and these would yield other subsets of inverse optimal cost functions of $\Gbf$. In the following section we will employ a convexifcation using a generalization of mutual information via $f$-divergence to obtain a larger class of 
inverse optimal cost functions for $\Gbf$.

\subsection{Inverse optimal control using $f$-divergence}
Let $f$ be a convex function such that $ f(1)=0 $.  The $f$-divergence of a pair of distributions $P,Q$ supported on the same finite set is given by $D_f(P||Q) = \sum_{x} Q(x) f\left( \frac{P(x)}{Q(x)} \right)$. One can define a corresponding $f$-mutual information between random variables $X,Y$ as:
\begin{align}
I_f(X;Y) &= \sum_{x,y} p(x,y) f\left(\frac{p(x)p(y)}{p(x,y)}\right)\non \\
 &= D_f(p_X(\cdot)p_Y(\cdot)||p(\cdot,\cdot)),
\end{align}
where $p(\cdot,\cdot)$ is the joint distribution of the random variables $X,Y$ and $p(x)$ (and $ p_X $), $p(y)$ (and $ p_Y $) are marginals of $X$ and 
$Y$. In this section we will derive inverse optimal cost functions for problem $\Gbf$ using a convexification that employs the $f$-mutual information. The K-L divergence is given by $D(P||Q) = \sum_{x} P(x) \log \frac{P(x)}{Q(x)}= D_{-\log }(Q||P),$ and therefore, mutual information is $I(X;Y) = I_{-\log } (X;Y)$. 

$f$-divergence enjoys many of the properties that K-L divergence has. Fundamental to this, is the following analogue of 
the $\log$-sum inequality, called the $f$-sum inequality~\cite{csiszar2004information}.
\begin{lemma}
Let $f$ be convex and let $a_i,b_i, i=1,\hdots,n$ be nonnegative. Then 
\[\sum_i b_i f\left(\frac{a_i}{b_i}\right) \geq \left(\sum_i b_i\right)f\left( \frac{\sum_i a_i}{\sum_i b_i}\right).\]
\end{lemma}

Furthermore, the $f$-mutual information satisfies the data processing inequality (for clarity, we call this the $f$-DPI). Since we were unable to locate a reference for this, we have included the proof for easy reference.
\begin{lemma}\label{lem:fdpi}
Suppose $X,Y,Z$ are discrete random variables such that $X \rightarrow Y \rightarrow Z$. Then, $I_f(X;Y) \geq I_f(X;Z)$.
\end{lemma}
\begin{IEEEproof}
On the one hand, we have
\begin{align*}
 I_f(X;(Y,Z)) &= \sum_{x,y,z} p(x,y,z) f\left(\frac{p(x)p(y,z)}{p(x,y,z)}\right),  \\
              &\buildrel{(a)}\over= \sum_{x,y,z} p(x,y,z)f\left(\frac{p(x)p(y,z)}{p(x|y)p(y,z)}\right),   \\
              &= \sum_{x,y,z} p(x,y,z) f\left(\frac{p(x)}{p(x|y)}\right)
         = I_f(X;Y).
\end{align*}
where in $(a)$ we have used that $X\rightarrow Y\rightarrow Z.$
On the other hand,
\begin{align*}
 I_f(X;&(Y,Z)) = \sum_{x,z} \sum_y p(x,y,z) f\left(\frac{p(x)p(y,z)}{p(x,y,z)}\right),  \\
              &\buildrel{(b)}\over\geq \sum_{x,z} \left( \sum_y p(x,y,z)\right)  f\left(\frac{\sum_y p(x)p(y,z)}{\sum_y p(x,y,z)}\right),   \\
              &= \sum_{x,z} p(x,z) f\left( \frac{p(x)p(z)}{p(x,z)} \right) 
              = I_f(X;Z),
\end{align*}
where in $(b)$ we have used the $f$-sum inequality. Combining these results we get the data processing inequality.
\end{IEEEproof}
Similarly if $X\rightarrow Y \rightarrow Z$, $I_f(Y;Z) \geq I_f(X;Z)$ and if $W \rightarrow X \rightarrow Y \rightarrow Z,$ then $I_f(X;Y) \geq I_f(W;Z)$. 

Just like the K-L divergence, $D_f(P||Q)$ is convex in $(P,Q)$~\cite{csiszar2004information}.
From this, we get the convexity of $I_f(X;Y)$ in the kernel $ p(y|x) $ for a fixed distribution $ p(x) $ of $X$.
\begin{lemma}
Let $ p(x,y) $ denote the joint distribution of random variables $ X,Y $ and $ p(x) $ denote the marginal of $ X. $ 
For fixed $p(x)$, $I_f(X;Y)$ is a convex function of $p(y|x)$. 
\end{lemma}
\begin{IEEEproof}
Let $ p_1,p_2 $ denote two distributions such that $ p_1(x)=p_2(x)=p(x) $. Let 
$ p_1(y|x), p_2(y|x) $ be the corresponding kernels and let $ \alpha \in [0,1]. $
Denote $p(y|x) = \alpha p_1(y|x) + (1-\alpha) p_2(y|x)$ and $ p(x,y) = p(y|x)p(x) $. It follows that $ p(x,y) = \alpha p_1(x,y) + (1-\alpha)p_2(x,y)$ and
$p(x)p(y) = p(x)\sum_{x'} p(y|x') p(x') = \alpha p(x)p_1(y) + (1-\alpha) p(x)p_2(y)$. Let $ I_f(p) $ and  $I_f(p_i) $ denote the mutual information of $ X,Y $  under joint distributions $ p(x,y) $ and $ p_i(x,y). $ We have
\begin{align*}
 I_f(p) &= D_f(p(x)p(y) || p(x,y)) \\
 &\buildrel{(c)}\over\leq \alpha D_f(p(x)p_1(y)||p_1(x,y)) \\ 
&\qquad + (1-\alpha)D_f(p(x)p_2(y)||p_2(x,y))\\
&= \alpha I_f(p_1) + (1-\alpha)I_f(p_2),
\end{align*}
where (c) is thanks to the convexity of $ D_f(\cdot||\cdot). $ This completes the proof.
\end{IEEEproof}

The convexity of the set $\Qscrdpi$ is due to the fact that the mutual information $I(X;Y)$ is convex in $p(y|x)$ for fixed $p(x)$ \textit{and} it is concave in $p(x)$ for fixed $p(y|x)$. The former property is shared by $f$-mutual information. 
For mutual information, the latter property is enabled by the additive decomposition of mutual information $I(X;Y) = H(Y) - H(Y|X)$. The additivity does not, in general, 
hold for $f$-mutual information and therefore it is not immediately clear if $I_f(X;Y)$ is concave in $p(x)$ for given $p(y|x)$. In section, we will restrict ourselves to those convex functions $f$ for which $f$-mutual information has this property.
\begin{definition}
A function $f$ is said to have the {\em saddle property} if $f$ is convex and if the $f$-mutual information $I_f(X;Y)$ is a concave function of $p(x)$ for any fixed $p(y|x)$. 
\end{definition}
\begin{remarkc}
$f(x) \equiv -\log(x)$ has the saddle property, so the set of such functions is not empty. Furthermore, this is not a vacuous definition $f(x) \equiv c -dx$ for constants $c,d$ also has the saddle property.
\end{remarkc}

Define the set
\begin{align}
\Qscrfdpi = \{Q \in \Pscr(\Zscr) | Q_S=P_S, Q_{Y|X}=P_{Y|X},\non \\
 I_f(Q_{S,\Shat}) \leq I_f(Q_{X,Y})\} 
\end{align}
By Lemma \ref{lem:fdpi}, $\Qscr \subset \Qscrfdpi$. Furthermore, we have the following lemma.
\begin{lemma}
Let $f$ have the saddle property. Then $\Qscrfdpi$ is convex.
\end{lemma}
\begin{IEEEproof}
Follows as in Lemma \ref{lem:convexg}.
\end{IEEEproof}
Consider the following relaxation of $ \Nbf, $ denoted $ f\td\Cbf_\Nbf. $
\begin{equation}
 \min_{Q \in \Qscrfdpi} \langle \kappa, Q \rangle.\tag{$f\td\Cbf_\Nbf$}
\end{equation} 
It follows that if $f$ has the saddle property the relaxation $f\td\Cbf_\Nbf$ is a convex optimization problem.
$f\td\Cbf_\Nbf$ can be written explicitly as a mathematical program in the variable $ Q \in \Qscrfdpi$ and additional variables $a =Q_{\shat|S},b=Q_X$, in the same way as $ \Cbf_\Nbf $ was written. 

Now consider problem $\Gbf$, \ie take $\kappa(z) = \delta(s,\shat)+\rho(x)$. Then following the same line of arguments made for $ \Cbf_\Gbf, $ the objective of the relaxation $ f\td\Cbf_\Gbf $
can be written in terms of $a,b$ alone and $Q$ can be eliminated.
$$ \problemsmall{$\Cbf_{\Gbf}$}
        {a,b}
        {\sum_{s,\shat}\delta(s,\shat)a(\shat|s)P_S(s) + \sum_{x}\rho(x)b(x)}
                                 {\hspace{-1.1cm}\begin{array}{r@{\, }c@{\, }l@{}c@{}}
\sum_{\shat}a(\shat|s)    & = & 1,  \forall s,\quad :\hspace{-1mm}\mu^a(s)P_S(s), & \sum_{x}b(x)             =  1, :\hspace{-1mm}\mu^b,\\ 
I_f(aP_S)                  & \leq  & I_f(P_{Y|X}b), \qquad :\lambda\\  
 a(\shat|s)  \geq &0,& \multicolumn{2}{l}{\hspace{-2mm}  \forall s,\shat\ \ :\hspace{-1mm}\nu^a(\shat|s)P_S(s),  \ \  b(x)     \hspace{-1mm}\geq  0,\, \forall x,\hspace{-.5mm} :\hspace{-1mm}\nu^b(x)}
        \end{array}} $$
Let $\mu^a,\mu^b, \lambda,\nu^a,\nu^b$ be the Lagrange multipliers of the constraints of problem $\Cbf_\Gbf$.  
$a,b$ are optimal for $\Cbf_\Gbf$ if and only if together with the Lagrange multipliers, they satisfy,
\begin{align}
 \delta(s,\shat)P_S(s)  =& -\lambda  \frac{\partial I_f(aP_S)}{\partial a(\shat|s)}+ P_S(s)\mu^a(s) + P_S(s)\nu^a(\shat|s) \tag{KKT$_{f\td\Cbf_\Gbf}$}\label{eq:kktfgbf}\\ 
 \rho(x)             =& \lambda  \frac{\partial I_f(P_{Y|X}b)}{\partial b(x)}+ \mu^b + \nu^b(x) \non \\ 
0\leq \lambda       \perp& I_f(P_{Y|X}b)-I_f(aP_S) \geq 0 \non\\
\textstyle \sum_{\shat}a(\shat|s) &= 1,  \left\langle a(\cdot|s), \nu^a(\cdot|s)\right\rangle    = 0, \nu^a(\shat|s)\geq 0, \non  \\
\textstyle \sum_{x}b(x)             =&  1, \left\langle b,\nu^b\right\rangle  =0, \nu^b(\cdot) \geq 0, b(\cdot) \geq 0,  a(\shat |s)\geq 0,\non
\end{align}
$\forall z\in \Zscr$. 
\begin{theorem}
Let $Q \in \Qscr$  be a candidate solution of $\Gbf$. Denote $a(\shat|s) = Q_{\Shat|S}(\shat|s), b(x) = Q_X(x)$. 
Then for this $Q$, the following functions $\delta,\rho$ are inverse optimal for problem $\Gbf$: 
\begin{align*}
 \delta(s,\shat)  =& -\frac{1}{P_S(s)}\lambda  \frac{\partial I_f(aP_S)}{\partial a(\shat|s)}+ \mu^a(s) + \nu^a(\shat|s) \\ 
 \rho(x)             =& \lambda  \frac{\partial I_f(P_{Y|X}b)}{\partial b(x)}+ \mu^b + \nu^b(x) \non 
\end{align*}
where $f$ is any continuously differentiable function having the saddle property, $\lambda \geq 0$ is a scalar satisfying
$0\leq \lambda       \perp I_f(P_{Y|X}b)-I_f(aP_S) \geq 0$
and the other parameters are constants or functions satisfying, $\mu^b \in \Real, \nu^b: \Xscr \rightarrow [0,\infty), \mu^a : \Sscr \rightarrow \Real, \nu^a : \Sscr \times \Sscrhat \rightarrow [0,\infty), \langle \nu^b, b \rangle =0, \langle \nu^a(\cdot|s), a(\cdot|s) \rangle =0$ for all $s$.
\end{theorem}
\begin{IEEEproof}
Under the given conditions, the system \eqref{eq:kktfgbf} is satisfied by $Q$. Furthermore, since $ f $ has the saddle property, $ f\td\Cbf_\Gbf $ is a convex optimization problem. This implies that $\delta$ and $\rho$ are inverse optimal for $Q$ for problem $f\td\Cbf_\Gbf$. As a consequence, they are inverse optimal for $\Gbf$. 
\end{IEEEproof}

\begin{remarkc}
Notice that the saddle property is not required to claim $ \Qscr \subset \Qscrfdpi.$ It is only required to claim that the $ Q $ that satisfies \eqref{eq:kktfgbf} is indeed optimal for $ f\td\Cbf_\Gbf, $ for which we need the convexity of $ f\td\Cbf_\Gbf. $
\end{remarkc}

\begin{remarkc}
$ \Qscrfdpi $ is a convex set that contains $ \Qscr$. Depending on the choice of $ f $ this relaxation may or may not be tight enough.
The tightest relaxation of $ \Qscr $ obtainable is the \textit{convex hull} of $ \Qscr. $ Characterizing this set remains a challenge.
\end{remarkc}
\section{CONVEX RELAXATION FOR THE WITSENHAUSEN PROBLEM} \label{sec:witsen}
We now study the relaxation with $\Qscrdpi$ of the Witsenhausen problem, \ie, $\Nbf$ with $\kappa(z) \equiv \zeta(x,\shat) + \tau(x,s)$. Consider the problems,\vspace{-30pt}
\begin{multicols}{2}
\begin{equation}
\minimize{Q \in \Qscr} \left\langle \zeta + \tau, Q \right\rangle, \tag{$\Wbf$}
\end{equation} \break
\begin{equation}
\minimize{Q \in \Qscrdpi} \left\langle \zeta + \tau, Q \right\rangle. \tag{$\Cbf_\Wbf$}
\end{equation}
\end{multicols} \vspace{-5pt}
\noindent The KKT conditions of $\Cbf_\Wbf$ are a special case of \eqref{eq:kktcbfnbf}. In Remark \ref{rem:kappastr}, we alluded to the fact that solving \eqref{eq:nu} is harder for this problem. Although we cannot solve $\Cbf_\Wbf$, we can characterize its optimal value in terms of the Langrage multiplier $\mu^a(\cdot)$. 
\begin{theorem}
Suppose $\Cbf_\Wbf$ has a solution $Q=Q^*$ and suppose system \eqref{eq:kktcbfnbf} characterizes the solution of $\Cbf_\Wbf$. Then the optimal value of $\Cbf_\Wbf$ is given by $-\left(\Ebb[\mu^{a*}(S)] + \lambda^*+ \mu^{b*}\right)$, where $\mu^{a*}(\cdot),\lambda^*,\mu^{b*}$ are the values of the (scaled) Langrange multipliers $\mu^a(\cdot),\lambda,\mu^b$ that satisfy \eqref{eq:kktcbfnbf} for $Q^*$.
\end{theorem}
\begin{IEEEproof}
In equation \eqref{eq:nu} of \eqref{eq:kktcbfnbf} take expectation with $Q^*$ to get
\begin{align*}
\Ebb[\zeta(X,\Shat) +\tau(X,S)&+\mu^{a*}(S)+ \lambda^* + \mu^{b*}] \\ &\qquad + \lambda (I(aP_S) - I(P_{Y|X}b))=0,
\end{align*}
where the expectation is over $Q^*$. Complementarity slackness for the DPI gives $\Ebb[\zeta(X,\Shat) + \tau(X,S)] = -\left(\Ebb[\mu^{a*}(S)]+\lambda^*+\mu^{b*}\right).$
\end{IEEEproof}
 $-\left(\Ebb[\mu^{a*}(S)]+\lambda^*+\mu^{b*}\right)$ is a lower bound on $\OPT(\Wbf)$ too. Notice that although $\mu^{a*}$ itself depends on 
$Q^*$, the expectation therein is over the \textit{given} distribution $P_S$. Finally note that the result applies for any cost structure and not only for $\Wbf$. 

We now show that it is not possible to solve \eqref{eq:kktcbfnbf} for $\Cbf_\Wbf$ with $\nu(\cdot)\equiv 0$ using linear controllers.
\begin{theorem} \label{thm:witsen}
 Consider the convex relaxation of the Witsenhausen problem, \ie, in $\Cbf_\Wbf$ suppose $S \sim \Nscr(0,\sigma_0^2)$, $Y= X+w$, $w \sim \Nscr(0,\sigma_w^2)$ and independent of $S$, $\zeta(x,\shat)\equiv (x-\shat)^2, \tau(x,s) \equiv (x-s)^2$. There exist no constants $\gamma_0,\gamma_1$ such that the PDF $(S,\gamma_0S,Y,\gamma_1Y)$ solves \eqref{eq:kktcbfnbf} with $\nu(\cdot)\equiv 0$.
\end{theorem}
\begin{IEEEproof}
Assume the contrary, \ie assume \eqref{eq:kktcbfnbf} holds for this PDF. The conditional distribution $a(\shat|s)$ induced by $\gamma_0,\gamma_1$ is positive for all $s,\shat$. Since $\langle a(\cdot|s),\nu^a(\cdot|s)\rangle=0$ for all $s$, $\nu^a(\shat|s) =0$ for almost every $\shat,s$. Likewise $\nu^b(x)=0$ for almost every $x$. Now in \eqref{eq:nu}, put $\nu(\cdot) =0$ and take expectation with $P_{Y|X}(\cdot|x)$ to get that for almost every $z$, 
$(x-\shat)^2  +(x-s)^2  + f_1(x) + f_2(s,\shat) + \mu^a(s)+c=0,$ where  $f_1(\cdot),f_2(\cdot)$ are quadratic polynomials,  $\mu^a(\cdot)$ is a function on $\Sscr  $ and $ c $ is a constant. This implies that $(x-\shat)^2 + f_2(s,\shat)$ is independent of  $\shat$ for arbitrary $x,s$ except from a set of measure zero. Clearly, this is not possible for a quadratic $f_2$.
\end{IEEEproof}
\begin{remarkc}[Linear controllers for the Witsenhausen problem:]
In a sense, the above theorem was only a sanity check. We already know from~\cite{witsenhausen_counterexample_1968} that linear controllers are not optimal for the Witsenhausen problem; had we found $\gamma_0,\gamma_1$ that solved $\Cbf_\Wbf$, 
they would automatically be optimal for the Witsenhausen problem. 
\end{remarkc}

This indicates that if at all a solution of the Witsenhausen problem is to be found through $\Cbf_\Wbf$, finding it will require a sophisticated attempt at solving \eqref{eq:kktcbfnbf}. Or, it may require a different convexification.
%

\section{CONCLUSIONS}\label{sec:conclusions}
This paper has presented an optimization-based approach and a general framework for stochastic control problems with nonclassical information structure that employed the idea of convex relaxation. We recovered solutions obtained by Bansal and Ba\c{s}ar through this approach. In our formulation, these stochastic control problems are nonconvex optimization problems and the information-theoretic data processing inequality appears as an artifice of convexifying them. We gave insights into the relation of the cost structure of the problem and the structure of the convexification for inverse optimal control. Using certain $ f $-divergences we obtained a wider set of inverse optimal cost functions than were obtainable through information theoretic methods.

\section{ACKNOWLEDGMENTS}
The authors would like to thank Profs Maxim Raginsky and P. R. Kumar for their inputs. Ankur's work was supported in part by the INSPIRE Faculty award of the Department of Science and Technology, Government of India.  The work of  T. P. Coleman was supported in part by the NSF Science \& Technology Center Grant CCF-0939370, in part by NSF Grant CCF-1065022, and in part by Systems on Nanoscale Information fabriCs (SONIC), one of the six SRC STARnet Centers, sponsored by MARCO and DARPA.







\bibliographystyle{IEEEtran}
\bibliography{ref_2,ref}
\end{document}

%% file: macros-ieee.tex
\usepackage{amsmath, amssymb, bbm, xspace}
\usepackage{epsfig}
\usepackage{longtable}
\usepackage{color}
\usepackage{mathrsfs}
\usepackage{subfigure}

\def\sthat{\ {\rm such\ that}\ }
\usepackage{courier}

\newtheorem{theorem}{Theorem}[section]

\newtheorem{lemma}{Lemma}[section]

\newtheorem{definition}{Definition}[section]
\def\bkE{{\rm I\kern-.17em E}}
\def\bk1{{\rm 1\kern-.17em l}}
\def\bkD{{\rm I\kern-.17em D}}
\def\bkR{{\rm I\kern-.17em R}}
\def\bkP{{\rm I\kern-.17em P}}

\def\bkZ{{\bf{Z}}}

\def\bkE{{\rm I\kern-.17em E}}
\def\bk1{{\rm 1\kern-.17em l}}
\def\bkD{{\rm I\kern-.17em D}}
\def\bkR{{\rm I\kern-.17em R}}
\def\bkP{{\rm I\kern-.17em P}}

\def\bkZ{{\bf{Z}}}
\def\b12{(\beta_1,\beta_2)}

\newcounter{example}
\renewcommand{\theexample}{\thesection.\arabic{example}}

\newcounter{remark}
\renewcommand{\theremark}{\thesection.\arabic{remark}}
\newenvironment{remarkc}[1][]{\refstepcounter{remark}
\noindent{\itshape Remark~\theremark. #1} \rmfamily}{\hspace*{\fill}~$\square$\vspace{0pt}}

\def\Xscr{\mathcal{X}}
\def\Yscr{\mathcal{Y}}
\def\indep{\perp\!\!\!\!\perp}
\def\Ebb{\mathbb{E}}
\newlength{\noteWidth}
\long\def\notes#1{\ifinner
{\tiny #1}
\else
\marginpar{\parbox[t]{\noteWidth}{\raggedright\tiny #1}}
\fi\typeout{#1}}

 \def\notes#1{\typeout{read notes: #1}} 

\newcommand{\I}[1]{\mathbb{I}_{\{#1\}}}

\newcommand{\ie}{i.e.\@\xspace} 
\newcommand{\eg}{e.g.\@\xspace} 
\newcommand{\etal}{et al.\@\xspace} 

\newcommand{\Real}{\ensuremath{\mathbbm{R}}}
\newcommand{\inv}{^{-1}}

\newcommand{\minimize}[1]{\displaystyle\minim_{#1}}
\newcommand{\minim}{\mathop{\hbox{\rm min}}}

\def\subject{\hbox{\rm s.t.}}

\def\OPT{{\rm OPT}}

\def\reff{^{\rm ref}}

\def\Gbf{{\bf G}}
\def\Cbf{{\bf C}}
\def\Bbf{{\bf B}}

\def\Wbf{{\bf W}}

\def\Nbf{{\bf N}}

\def\sgn{\mathop{\hbox{\rm sgn}}}

\def\norm#1{\|#1\|}
\def\spose#1{\hbox to 0pt{#1\hss}}

\def\text #1{\hbox{\quad#1\quad}}

\def\nthinsp{\mskip -2   mu}

\def\Q{_{\scriptscriptstyle Q}}

\def\superstar{^{\raise 0.5pt\hbox{$\nthinsp *$}}}
\def\SUPERSTAR{^{\raise 0.5pt\hbox{$*$}}}

\def\lamstarT {\lambda^{\raise 0.5pt\hbox{$\nthinsp *$}T}}

\def\s{^*}

\def\del{\partial}

\def\Pscr{{\cal P}}
\def\Qscr{{\cal Q}}
\def\Sscr{{\cal S}}
\def\Sscrhat{\widehat{\mathcal{S}}}

\def\Nscr{{\cal N}}

\def\Zscr{{\cal Z}}

\def\shat{\widehat s}

\def\Shat{\widehat S}

\def\Nbf{{\bf N}}

\def\aur{\;\textrm{and}\;}

\def\del{\partial}
\def\non{\nonumber}

\let\forallnew\forall
\renewcommand{\forall}{\forallnew\ }
\let\forall\forallnew

\def\ds{\displaystyle}
		\def\bkE{{\rm I\kern-.17em E}}
		\def\bk1{{\rm 1\kern-.17em l}}
		\def\bkD{{\rm I\kern-.17em D}}
		\def\bkR{{\rm I\kern-.17em R}}
		\def\bkP{{\rm I\kern-.17em P}}
		\def\bkY{{\bf \kern-.17em Y}}
		\def\bkZ{{\bf \kern-.17em Z}}
		\def\bkC{{\bf  \kern-.17em C}}

{\begin{list}{}%
         {\setlength{\leftmargin}{#1}}%
         \item[]%
}
{\end{list}}

\def\Dsf{\mathsf{D}}

		\def\bsp{\begin{split}}
		\def\beq{\begin{eqnarray}}
		\def\bal{\begin{align*}}
		\def\bc{\begin{center}}
		\def\be{\begin{enumerate}}
		\def\bi{\begin{itemize}}
		\def\bs{\begin{small}}
		\def\bS{\begin{slide}}
		\def\ec{\end{center}}
		\def\ee{\end{enumerate}}
		\def\ei{\end{itemize}}
		\def\es{\end{small}}
		\def\eS{\end{slide}}
		\def\eeq{\end{eqnarray}}
		\def\eal{\end{align*}}
		\def\esp{\end{split}}
		\def\problemlarge#1#2#3#4{\fbox
		 {\begin{tabular*}{.88\textwidth}
			{@{}l@{\extracolsep{\fill}}l@{\extracolsep{6pt}}l@{\extracolsep{\fill}}c@{}}
				#1 & $\minimize{#2}$ & $#3$ & $ $ \\[5pt]
					 & $\subject\ $    & $#4$ & $ $
			\end{tabular*}}
			}

		\def\problemsmall#1#2#3#4{\fbox
		 {\begin{tabular*}{0.47\textwidth}
			{@{}l@{\extracolsep{\fill}}l@{\extracolsep{6pt}}l@{\extracolsep{\fill}}c@{}}
				#1 & $\minimize{#2}$ & $#3$ & $ $ \\[5pt]
					  $\subject\ $ &    & $#4$ & $ $
			\end{tabular*}}
			}

	\def\cp2problem#1#2#3#4{\fbox
		 {\begin{tabular*}{0.9\textwidth}
			{@{}l@{\extracolsep{\fill}}l@{\extracolsep{6pt}}l@{\extracolsep{\fill}}c@{}}
				#1 & & $#4 $ 
			\end{tabular*}}}

\newcommand{\pmat}[1]{\begin{pmatrix} #1 \end{pmatrix}}

		\def\bkE{{\rm I\kern-.17em E}}
		\def\bk1{{\rm 1\kern-.17em l}}
		\def\bkD{{\rm I\kern-.17em D}}
		\def\bkR{{\rm I\kern-.17em R}}
		\def\bkP{{\rm I\kern-.17em P}}
		
		\def\bkZ{{\bf{Z}}}

\newcommand {\beeq}[1]{\begin{equation}\label{#1}}
\newcommand {\eeeq}{\end{equation}}
\newcommand {\bea}{\begin{eqnarray}}
\newcommand {\eea}{\end{eqnarray}}

\def\texitem#1{\par\smallskip\noindent\hangindent 25pt
               \hbox to 25pt {\hss #1 ~}\ignorespaces}

